\documentclass[12pt,seceqn,secthm,amsthm]{elsart}%
\usepackage{amsfonts}
\usepackage{amsmath}
\usepackage{amssymb}
\usepackage[dvips]{graphicx}%
\def\openone{\hbox{\upshape \small1\kern-3.8pt\normalsize1}}
\numberwithin{figure}{section}
\newcounter{savecounter}
\renewcommand{\theenumi}{\roman{enumi}}

\begin{document}
%
\begin{frontmatter}%
%

\title{Measures in wavelet decompositions\protect\rule{0pt}{96pt}}%
%

\author{Palle E. T. Jorgensen\thanksref{label1}}\thanks
[label1]{This material is based upon work supported by the U.S.
National Science Foundation under Grant No.\ DMS-0139473 (FRG).}%

%

\address{Department of Mathematics,
The University of Iowa,
14 MacLean Hall,
Iowa City, IA 52242-1419,
U.S.A.}%
\ead{jorgen@math.uiowa.edu}%
\ead[url]{http://www.math.uiowa.edu/\symbol{126}jorgen}%
%

\begin{abstract}
In applications, choices of orthonormal bases in Hilbert space $\mathcal{H}$
may come about from the simultaneous diagonalization of some specific abelian
algebra of operators. This is the approach of quantum theory as suggested by
John von Neumann; but as it turns out, much more recent constructions of bases
in wavelet theory, and in dynamical systems, also fit into this scheme.
However, in these modern applications, the basis typically comes first, and
the abelian algebra might not even be made explicit. It was noticed recently
that there is a certain finite set of non-commuting operators $F_{i}$, first
introduced by engineers in signal processing, which helps to clarify this
connection, and at the same time throws light on decomposition possibilities
for wavelet packets used in pyramid algorithms. There are three interrelated
components to this: an orthonormal basis, an abelian algebra, and a
projection-valued measure. While the operators $F_{i}$ were originally
intended for quadrature mirror filters of signals, recent papers have shown
that they are ubiquitous in a variety of modern wavelet constructions, and in
particular in the selection of wavelet packets from libraries of bases. These
are constructions which make a selection of a basis with the best frequency
concentration in signal or data-compression problems. While the algebra
$\mathcal{A}$ generated by the $F_{i}$-system is non-abelian, and goes under
the name \textquotedblleft Cuntz algebra\textquotedblright\ in $C^{\ast}%
$-algebra theory, each of its representations contains a canonical maximal
abelian subalgebra, i.e., the subalgebra is some $C(X)$ for a Gelfand space
$X$. A given representation of $\mathcal{A}$, restricted to $C(X)$, naturally
induces a projection-valued measure on $X$, and each vector in $\mathcal{H}$
induces a scalar-valued measure on $X$. We develop this construction in the
general context with a view to wavelet applications, and we show that the
measures that had been studied earlier for a very restrictive class of $F_{i}%
$-systems (i.e., the Lemari\'{e}-Meyer quadrature mirror filters) in the
theory of wavelet packets are special cases of this. Moreover, we prove a
structure theorem for certain classes of induced scalar measures. In the
applications, $X$ may be the unit interval, or a Cantor set; or it may be an
affine fractal, or even one of the more general iteration limits involving
iterated function systems consisting of conformal maps.
\end{abstract}%
%

\begin{keyword}
Hilbert space, Cuntz algebra, completely positive map, creation operators,
wavelet packets, pyramid algorithm, product measures, orthogonality relations,
equivalence of measures, iterated function systems (IFS), scaling function,
multiresolution, subdivision scheme, singular measures, absolutely continuous
measures
\renewcommand{\MSC}{{\par\leavevmode\hbox{\it2000 MSC:\ }}}\MSC
42C40; 42A16; 43A65; 42A65%
\end{keyword}%
%

\end{frontmatter}%

\section{\label{Int}Introduction}

A popular approach to wavelet constructions is based on a so-called scaling
identity, or scaling equation. A solution to this equation is a function on
$\mathbb{R}^{d}$ for some $d$. The equation is related to a subdivision scheme
that is used in numerical analysis and in computer graphics. In that language,
it arises from a fixed scaling matrix, assumed expansive, a system of masking
coefficients, and a certain subdivision algorithm. An iteration of the scaling
produces a succession of subdivisions into smaller and smaller frequency
bands. In signal processing, the coefficients in the equation refer to
``frequency response''. There are various refinements, however, of this setup:
two such refinements are \emph{multi-wavelets} and \emph{singular systems}.

If the masking coefficients are turned into a generating function, called a
low-pass filter $m_{0}$, then the scaling identity takes a form which admits
solutions with an infinite product representation. Various regularity
assumptions are usually placed on the function $m_{0}$. The first requirement
is usually that the solution, i.e., the scaling function, is in $L^{2}%
(\mathbb{R}^{d})$, but other Hilbert spaces of functions on $\mathbb{R}^{d}$
are also considered. If the number of masking coefficients is finite, then
$m_{0}$ is a Fourier polynomial. (For the Daubechies wavelet, there are four
coefficients, and $d=1$.) Readers not familiar with wavelets are referred to
the classic \cite{Dau92} by Daubechies. More general families of
multiresolutions are studied in \cite{PaRi03}, \cite{Rus+92}, and
\cite{DHJLW02}. For recent applications of multiresolutions to physics, see
\cite{JoPa03}. In general, however, $m_{0}$ might be a fairly singular
function. In favorable cases, the associated infinite product will be the
Fourier transform of the scaling function. This function, sometimes called the
father function, is the starting point of most wavelet constructions, the
multiresolution schemes. The function $m_{0}$ is a function of one or more
frequency variables, and convergence of the associated infinite product
dictates requirements on $m_{0}$ for small frequencies, hence low-pass. The
term \textquotedblleft low-pass\textquotedblright\ suggests a filter which
lets low-frequency signals pass with high probability. A complete system, of
which $m_{0}$ is a part, and which is built from appropriately selected
frequency bands, offers an effective tool for wavelet analysis and for signal
processing. Such a system gives rise to operators $F_{i}$, and their duals
$F_{i}^{\ast}$, that are the starting point for a class of algorithms called
pyramid algorithms. They are basic to both signal processing and the analysis
of wavelet packets. (In operator theory, $F_{i}^{\ast}$ is usually denoted
$S_{i}$, and $S_{i}^{\ast}$ is set equal to $F_{i}$. The reason is that it is
the operator $F_{i}^{\ast}$ that is isometric.) In the more traditional
approaches, $m_{0}$ is a Fourier polynomial, or at least a Lipschitz-class
function on a suitable torus, and the low-pass signal analysis is then
relatively well understood. But a variety of applications, for example to
multi-wavelets, dictate filters $m_{0}$ that are no better than continuous, or
perhaps only measurable. Then the standard tools break down, and probabilistic
and operator theoretic methods are forced on us. This is the setting which is
the focus of the present paper.

Recent developments in wavelet analysis have brought together ideas
from engineering and from computational mathematics, as well as fundamentals
from representation theory. One of the aims of this paper is to stress the
interconnections, as opposed to one aspect of this in isolation. 

By now, the subject draws on ideas from a variety of directions. Of these
directions, we single out quadrature-mirror filters from signal/image
processing, see Figure \ref{FigPerfectReconstruction} below.
High-pass/low-pass signal-processing
algorithms have now been adopted by pure mathematicians, although they
historically first were intended for speech signals, see \cite{Jor03}.
Perhaps
unexpectedly, essentially the same quadrature relations were rediscovered in
operator algebra theory, and they are now used in relatively painless
constructions of varieties of wavelet bases. The connection to signal
processing is rarely stressed in the math literature. Yet, the flow of ideas
between signal processing and wavelet mathematics
is a success story that deserves to be told.
Without these recent synergistic trends, we would perhaps only
know isolated examples of wavelets. Thus, mathematicians have borrowed from
engineers; and the engineers may be happy to know that what they do is used
in mathematics.

Our new results in this paper include
Corollary \ref{CorProNew.1},
Proposition \ref{ProCom.2},
Theorem \ref{ThmFam.3}, and
Corollary \ref{CorCom.4}, 
covering both construction (algorithms) for
wavelets, and selection (statistics) of the ``best'' wavelets in explicitly
parametrized families. They concern a construction of measures which allows
the selection of the ``best'' wavelet from a library of wavelet bases
(decomposition theory).

It is well known that the quadrature mirror filters which are used in subband
constructions of signal processing are also the building blocks for wavelets
and for wavelet packets; see, e.g., \cite{CMW92b} and \cite{Jor03}. The reader
may find good accounts of recent results on wavelet packets in the papers
\cite{DHJLW02} and \cite{Rus+92}. The scaling function for the wavelets, and
the wavelet packet functions arise from pyramid algorithms which are built
directly from the quadrature mirror filters. While the wavelet functions live
in spaces of functions on $\mathbb{R}$, typically $L^{2}\left(  \mathbb{R}%
\right)  $, the signals may be analyzed in the sequence space $\ell^{2}\left(
\mathbb{Z}\right)  $, or equivalently $L^{2}\left(  \mathbb{T}\right)  $,
where $\mathbb{T}=\mathbb{R}/2\pi\mathbb{Z}$. As is well known, the
isomorphism $L^{2}\left(  \mathbb{T}\right)  \cong\ell^{2}\left(
\mathbb{Z}\right)  $ is given by the transform of Fourier series. Then there
is an operator which maps $\ell^{2}\left(  \mathbb{Z}\right)  $ onto some
resolution subspace in $L^{2}\left(  \mathbb{R}\right)  $ and intertwines the
analysis of the signals in $\ell^{2}$ with the transformations acting on the
wavelet functions. In the simplest case, there is a function $\varphi\in
L^{2}\left(  \mathbb{R}\right)  $, called the scaling function, which sets up
the operator from $\ell^{2}$ to $L^{2}\left(  \mathbb{R}\right)  $: If
$\xi=\left(  \xi_{k}\right)  _{k\in\mathbb{Z}}\in\ell^{2}$, set
\begin{equation}
\left(  W_{\varphi}\xi\right)  \left(  x\right)  =\sum_{k\in\mathbb{Z}}\xi
_{k}\varphi\left(  x-k\right)  . \label{eqInt.1}%
\end{equation}
A subband filter is given by a sequence $\left(  a_{k}\right)  _{k\in
\mathbb{Z}}$ of \emph{frequency response coefficients}. They define an
operator $S_{0}$ on $\ell^{2}$ as follows:%
\begin{equation}
\left(  S_{0}\xi\right)  _{n}=\sum_{k}a_{n-2k}\xi_{k} \,, \label{eqInt.2}%
\end{equation}
and it is denoted [filter]$\bigcirc\llap{\small$\uparrow\mkern5.25mu$}$, i.e.,
it is a composition of the two operations, with $\bigcirc
\llap{\small$\uparrow\mkern5.25mu$}$ being the symbol for up-sampling; see
\cite{StNg96} and \cite{BrJo02b} for details. A function $\varphi$ on
$\mathbb{R}$ is said to satisfy a scaling identity with \emph{masking
coefficients} $\left(  a_{k}\right)  _{k\in\mathbb{Z}}$ if%
\begin{equation}
\varphi\left(  x\right)  =\sqrt{2}\sum_{k\in\mathbb{Z}}a_{k}\varphi\left(
2x-k\right)  . \label{eqInt.3}%
\end{equation}
The following lemma makes the connection between the discrete analysis of
$\ell^{2}$ and the wavelet analysis on $\mathbb{R}$.

The issue of smoothness properties of the possible scaling functions $\varphi
$, and the corresponding wavelets, is an important one. It is studied in a
number of papers, for example in \cite{CoDa96}, and the reader will find more
in \cite{BrJo02b}.

\begin{lem}
\label{LemInt.1}Suppose the sequence $\left(  a_{k}\right)  _{k\in\mathbb{Z}}$
is such that the function%
\begin{equation}
m_{0}\left(  z\right)  =\sum_{k}a_{k}z^{k} \label{eqInt.4}%
\end{equation}
is in $L^{\infty}\left(  \mathbb{T}\right)  $. Let $S_{0}$ be the
corresponding bounded operator on $\ell^{2}$. Let $\varphi\in L^{2}\left(
\mathbb{R}\right)  $, and let $W_{\varphi}$ be the corresponding operator
\textup{(\ref{eqInt.1})}. Then $\varphi$ satisfies the scaling identity
\textup{(\ref{eqInt.3})} if and only if%
\begin{equation}
W_{\varphi}S_{0}\xi=\frac{1}{\sqrt{2}}\left(  W_{\varphi}\xi\right)  \left(
\frac{x}{2}\right)  . \label{eqInt.5}%
\end{equation}

\end{lem}

In other words, $W_{\varphi}$ intertwines $S_{0}$ with the dyadic scaling
operator on $L^{2}\left(  \mathbb{R}\right)  $. We shall introduce
\begin{equation}
\left(  Uf\right)  \left(  x\right)  =\frac{1}{\sqrt{2}}f\left(  \frac{x}%
{2}\right)  \label{eqUnitaryScaling}%
\end{equation}
for the unitary scaling operator on $L^{2}\left(  \mathbb{R}\right)  $, and
\textup{(\ref{eqInt.5})} takes the form%
\begin{equation}
W_{\varphi}S_{0}=UW_{\varphi}. \label{eqInt.6}%
\end{equation}

\begin{pf}
The proof is straightforward, and we refer to \cite{BrJo02b} or \cite{JoKr03}
for details. \qed
\end{pf}

The quadrature conditions on the filter $\left(  a_{k}\right)  $ may be stated
as%
\begin{equation}
\sum_{k}\bar{a}_{k}a_{k+2l}=\delta_{0,l},\qquad l\in\mathbb{Z}.
\label{eqInt.7}%
\end{equation}
If
\begin{equation}
m_{0}\left(  z\right)  =\sum_{k}a_{k}z^{k}\text{\quad and\quad}m_{1}\left(
z\right)  =z\,\overline{m_{0}\left(  -z\right)  }\,,\qquad z\in\mathbb{T}%
\text{,} \label{eqInt.8}%
\end{equation}
then the two operators $S_{0}$ and $S_{1}$ given by
\begin{multline}
\left(  S_{i}f\right)  \left(  z\right)  =m_{i}\left(  z\right)  f\left(
z^{2}\right)  ,\\
f\in L^{2}\left(  \mathbb{T}\right)  ,\;z\in\mathbb{T}=\left\{  \,
z\in\mathbb{C}\mid\left\vert z\right\vert =1\,\right\}  ,\;i=0,1, \label{eqSi}%
\end{multline}
define isometries on $L^{2}\left(  \mathbb{T}\right)  $, and they satisfy the
relations%
\begin{align}
\sum_{i}S_{i}S_{i}^{\ast}  &
=\hbox{\upshape \small1\kern-3.8pt\normalsize1}_{L^{2}\left(  \mathbb{T}%
\right)  },\label{eqInt.9}\\
S_{i}^{\ast}S_{j}  &  =\delta_{i,j}%
\hbox{\upshape \small1\kern-3.8pt\normalsize1}_{L^{2}\left(  \mathbb{T}%
\right)  }, \label{eqInt.10}%
\end{align}
where $\hbox{\upshape \small1\kern-3.8pt\normalsize1}_{L^{2}\left(
\mathbb{T}\right)  }$ denotes the identity operator in the Hilbert space
$L^{2}\left(  \mathbb{T}\right)  $, and $\mathbb{T}$ is equipped with the
usual Haar measure. Equivalently, $L^{2}\left(  \mathbb{T}\right)  $ is viewed
as a space of $2\pi$-periodic functions, and the measure on $\mathbb{T}$ is
then $\left(  2\pi\right)  ^{-1}\,d\theta$. The relations (\ref{eqInt.9}%
)--(\ref{eqInt.10}) are called the Cuntz relations, see Section \ref{Pro}
below, but they also reflect the realization of the diagram in Figure
\ref{FigPerfectReconstruction}, from signal processing.

\begin{figure}[ptb]
\begin{center}
\setlength{\unitlength}{2pc} \begin{picture}(16.2,9)(0,-4.5)
\put(0,0){\vector(1,0){3}}
\put(0,0){\makebox(3,1){\textsc{signal in}}}
\put(0,-1){\makebox(3,1){$\xi$}}
\put(3,0){\line(3,5){1.5}}
\put(4.5,2.5){\line(1,0){1.25}}
\put(3,0){\line(3,-5){1.5}}
\put(4.5,-2.5){\line(1,0){1.25}}
\put(3,3.5){\line(1,0){4.25}}
\put(3,3.5){\line(0,-1){0.25}}
\put(7.25,3.5){\line(0,-1){0.25}}
\put(3,3){\makebox(4.25,0.5){\textsc{analysis}}}
\put(3,3.5){\makebox(4.25,1){$S_0^*$}}
\put(5.5,2){\makebox(1,1){$\bigcirc \llap{\small$\downarrow\mkern5.25mu$}$}}
\put(2.75,1){\makebox(1,1)[r]{low-pass filter}}
\put(2.75,-2){\makebox(1,1)[r]{high-pass filter}}
\put(5.5,-3){\makebox(1,1){$\bigcirc \llap{\small$\downarrow\mkern5.25mu$}$}}
\put(3,-3.5){\line(1,0){4.25}}
\put(3,-3.5){\line(0,1){0.25}}
\put(7.25,-3.5){\line(0,1){0.25}}
\put(3,-4.5){\makebox(4.25,1){$S_1^*$}}
\put(6.25,2.5){\line(1,0){3}}
\put(6.25,-2.5){\line(1,0){3}}
\put(5.5,-0.5){\makebox(1,1){down-sampling}}
\put(13.2,0){\vector(1,0){3}}
\put(13.2,0){\makebox(3,1){\textsc{signal out}}}
\put(13.2,-1){\makebox(3,1){$\xi$}}
\put(11,2.5){\line(3,-5){1.2}}
\put(9.75,2.5){\line(1,0){1.25}}
\put(11,-2.5){\line(3,5){1.2}}
\put(9.75,-2.5){\line(1,0){1.25}}
\put(9,2){\makebox(1,1){$\bigcirc \llap{\small$\uparrow\mkern5.25mu$}$}}
\put(11.75,1){\makebox(1,1)[l]{dual low-pass filter}}
\put(11.75,-2){\makebox(1,1)[l]{dual high-pass filter}}
\put(9,-3){\makebox(1,1){$\bigcirc \llap{\small$\uparrow\mkern5.25mu$}$}}
\put(9,-0.5){\makebox(1,1){up-sampling}}
\put(12.2,-0.5){\framebox(1,1){\huge$+$}}
\put(8.25,3.5){\line(1,0){4.25}}
\put(8.25,3.5){\line(0,-1){0.25}}
\put(12.5,3.5){\line(0,-1){0.25}}
\put(8.25,3){\makebox(4.25,0.5){\textsc{synthesis}}}
\put(8.25,3.5){\makebox(4.25,1){$S_0$}}
\put(8.25,-3.5){\line(1,0){4.25}}
\put(8.25,-3.5){\line(0,1){0.25}}
\put(12.5,-3.5){\line(0,1){0.25}}
\put(8.25,-4.5){\makebox(4.25,1){$S_1$}}
\end{picture}
\end{center}
\caption{Perfect reconstruction of signals}%
\label{FigPerfectReconstruction}%
\end{figure}
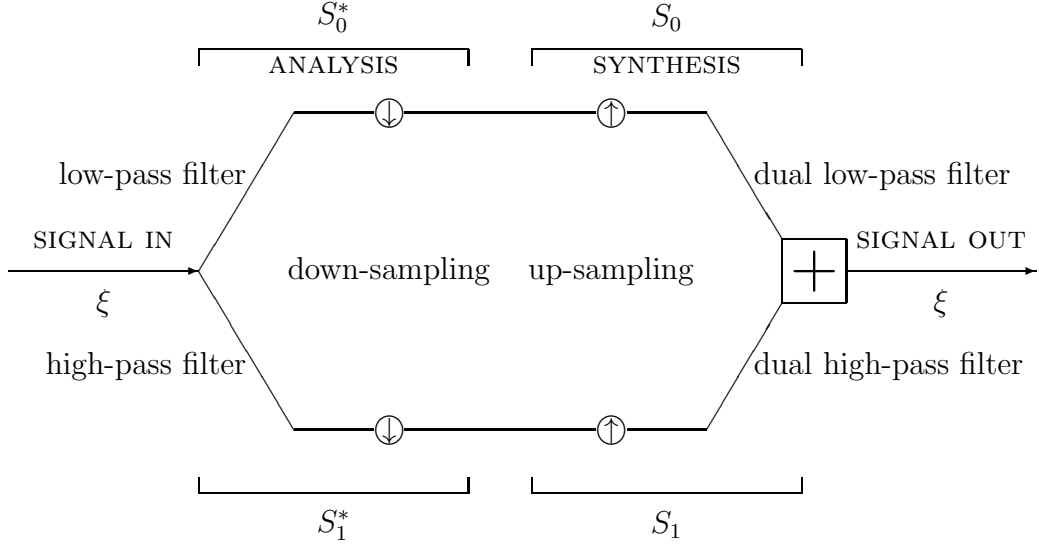

When the two operators and their dual adjoints act on sequences, then
(\ref{eqInt.9}) takes the form%
\begin{equation}
S_{0}S_{0}^{\ast}\xi+S_{1}S_{1}^{\ast}\xi=\xi\label{eqInt.11}%
\end{equation}
and expresses perfect reconstruction of signals from the subbands.

In view of (\ref{eqInt.9})--(\ref{eqInt.10}) it is clear that the isometries
$S_{i}$ provide dyadic subdivisions of the Hilbert space $\mathcal{H}%
=L^{2}\left(  \mathbb{T}\right)  \cong\ell^{2}$. Specifically, for every
$k\in\mathbb{Z}_{+}$ the subspaces%
\begin{equation}
\mathcal{H}\left(  i_{1},i_{2},\dots,i_{k}\right)  :=S_{i_{1}}S_{i_{2}}\cdots
S_{i_{k}}\mathcal{H} \label{eqInt.12}%
\end{equation}
are mutually orthogonal, and%
\begin{equation}
\sideset{}{^{\smash{\oplus}}}{\sum}\limits_{i_{1},\dots,i_{k}}\mathcal{H}%
\left(  i_{1},i_{2},\dots,i_{k}\right)  =\mathcal{H}. \label{eqInt.13}%
\end{equation}
But if the index labels $\left(  i_{1},\dots,i_{k}\right)  $ are used in
assigning dyadic partitions, for example the intervals $\displaystyle\left[
\frac{i_{1}}{2}+\dots+\frac{i_{k}}{2^{k}},\frac{i_{1}}{2}+\dots+\frac{i_{k}%
}{2^{k}}+\frac{1}{2^{k}}\right)  $, then it can be shown that $\left(
i_{1},\dots,i_{k}\right)  \mapsto\mathcal{H}\left(  i_{1},\dots,i_{k}\right)
$ extends to a projection-valued measure $E$ on the unit interval $I$, defined
on the Borel subsets of $I$, specifically%
\begin{equation}
E\left(  \left[  \frac{i_{1}}{2}+\dots+\frac{i_{k}}{2^{k}},\frac{i_{1}}%
{2}+\dots+\frac{i_{k}}{2^{k}}+\frac{1}{2^{k}}\right)  \right)  =\mathcal{H}%
\left(  i_{1},\dots,i_{k}\right)  , \label{eqInt.14}%
\end{equation}
and it was shown in \cite{CMW92b} and \cite{Sal03} that this measure
determines the selection of bases of wavelet packets from some prescribed
library of bases. The libraries of bases in turn are determined by quadrature
mirror filters.

However, it is difficult to compute $E\left(  \,\cdot\,\right)  $ in general.
If $f\in\mathcal{H}$, $\left\Vert f\right\Vert =1$, then%
\begin{equation}
\mu_{f}\left(  \,\cdot\,\right)  :=\left\langle \, f\mid E\left(
\,\cdot\,\right)  f\,\right\rangle =\left\Vert E\left(  \,\cdot\,\right)
f\right\Vert ^{2} \label{eqInt.15}%
\end{equation}
is a probability measure on $I$, and it is easier to compute for special
classes of quadrature mirror filters; explicit results are given in
\cite{CMW92b} and \cite{Sal03} for the filters $m_{0}$, $m_{1}$ first
introduced by Y.~Meyer. But it is not known in general for which quadrature
mirror filters $m_{i}$, and for which $f\in\mathcal{H}$, the measure $\mu
_{f}\left(  \,\cdot\,\right)  =\left\Vert E\left(  \,\cdot\,\right)
f\right\Vert ^{2}$ is absolutely continuous with respect to Lebesgue measure
on $I$. Absolute continuity is desirable in the calculus of libraries of bases
formed from wavelet packets.

\begin{rem}
\label{Remark1.2}While the conditions we list in (\ref{eqInt.9})--(\ref{eqInt.10}) may seem
unnecessarily stringent, it is possible to use the methods in our paper on a wider class of operator systems $S_i$ than the ones which correspond to perfect reconstruction, as we define 
it by Figure \ref{FigPerfectReconstruction}. In fact, Arveson \cite{Arv04} has recently 
developed an elegant operator-theoretic approach to finite systems of operators 
$S_i$,  $i = 0, 1, \dots, n$, when it is only assumed that the operator system of $n + 1$ operators
forms a \emph{row-contraction}. By this we mean that each operator $S_i$ is defined in a 
Hilbert space $\mathcal{H}$, and the system satisfies  
the contractivity condition
\[
\left\Vert
\sum_{i=0}^{n}S_{i}f_{i}
\right\Vert^{2}
\le \,
\sum_{i=0}^{n}
\left\Vert
f_{i}
\right\Vert^{2}
\text{\qquad for all }
\left( f_{0},\dots,f_{n}\right) \in \bigoplus_{0}^{n}\mathcal{H},
\]
or equivalently
\[
\sum_{i=0}^{n}
\left\Vert
S_{i}^{\ast}f
\right\Vert^{2}
\le
\left\Vert
f
\right\Vert^{2}
\text{\qquad for all }
f\in\mathcal{H}.
\]
As stressed in papers by Ron and Shen, e.g., \cite{RoSh00}, such row-contractions arise from
conditions on systems of filter functions which are weaker than the ones we summarize in 
equations (\ref{eqInt.7})--(\ref{eqInt.8}) above. The corresponding function system in 
$L^{2}\left(  \mathbb{R}\right)  $ will then not be a wavelet system in the sense we discuss 
below. It will only have considerably weaker orthogonality properties than those which are 
customary for the standard wavelet bases, and the authors of \cite{RoSh00} refer to these 
systems as \emph{framelets}; see also our survey paper \cite{Jor04d}.
\end{rem}

\section{\label{Sub}Subdivisions}

Subdivisions serve as an effective tool in the theory of dynamical systems
\cite{Rue94}, in computations \cite{SHS+99}, and in approximation theory
\cite{RoSh00}; see also \cite{BrJo02b}. Moreover, they are part of many
wavelet constructions: see, e.g., \cite{SSZ99}. The simplest such is the
familiar representations of the fractions $0\leq x<1$ in base $2$. For
$k\in\mathbb{N}$ and $a_{1},a_{2},\dots,a_{k}\in\left\{  0,1\right\}  $, set%
\begin{equation}
J_{k}\left(  a\right)  :=\left[  \frac{a_{1}}{2}+\dots+\frac{a_{k}}{2^{k}%
},\frac{a_{1}}{2}+\dots+\frac{a_{k}}{2^{k}}+\frac{1}{2^{k}}\right)  .
\label{eqSub.1}%
\end{equation}
Each interval $J_{k}\left(  a\right)  $ is contained in some $J_{k-1}\left(
b\right)  $, and the length of $J_{k}\left(  a\right)  $ is $2^{-k}$ by
definition. Moreover, the symbols $\left(  a_{1},\dots,a_{k}\right)  $, for
$k$ finite, yield a one-to-one representation of the dyadic rational
fractions. Note that we are excluding those infinite strings which terminate
with an infinite tail of $1$'s, and an infinite tail of $0$'s may be omitted
in listing the bits $a_{1},a_{2},\dots,a_{k}$.

We will also need the analogous representation of fractions in base $N$ where
$N\in\mathbb{Z}_{+}$, $N\geq2$. In that case $a_{i}\in\left\{  0,1,\dots
,N-1\right\}  $, the left-hand endpoint of $J_{k}\left(  a\right)  $ is
$\displaystyle\frac{a_{1}}{N}+\dots+\frac{a_{k}}{N^{k}}$, and
$\operatorname*{length}\left(  J_{k}\left(  a\right)  \right)  =N^{-k}$.

More general partitions like this arise in the study of endomorphisms
$\sigma\colon X\rightarrow X$, where $X$ is a compact Hausdorff space, and
$\sigma$ is continuous and onto. If, for each $x\in X$, the cardinality of
$\sigma^{-1}\left(  x\right)  =\left\{  \, y\in X\mid\sigma\left(  y\right)
=x\,\right\}  $ is $N$, independently of $x$, then there are branches of the
inverse, i.e., maps
\begin{equation}
\sigma_{0},\sigma_{1},\dots,\sigma_{N-1}\colon X\longrightarrow X
\label{eqSub.2}%
\end{equation}
such that
\begin{equation}
\sigma\circ\sigma_{i}=1_{X}, \label{eqSub.2a}%
\end{equation}
or in other notation,%
\begin{equation}
\sigma\left(  \sigma_{i}\left(  x\right)  \right)  =x,\qquad x\in X,
\label{eqSub.2b}%
\end{equation}
for $0\leq i<N$. Naturally, it is of special interest if the sections
$\left\{  \sigma_{i}\right\}  _{0\leq i<N}$ may be chosen to be continuous, as
is the case in the study of complex iteration of rational maps; see, e.g.,
\cite{CaGa93}, \cite{YHK97}.

\begin{figure}[ptb]
\begin{center}
\setlength{\unitlength}{0.3bp} \begin{picture}(390,390)(-15,-15)
\put(0,1){\includegraphics[bb=0 0 360 360,width=108bp,height=108bp]{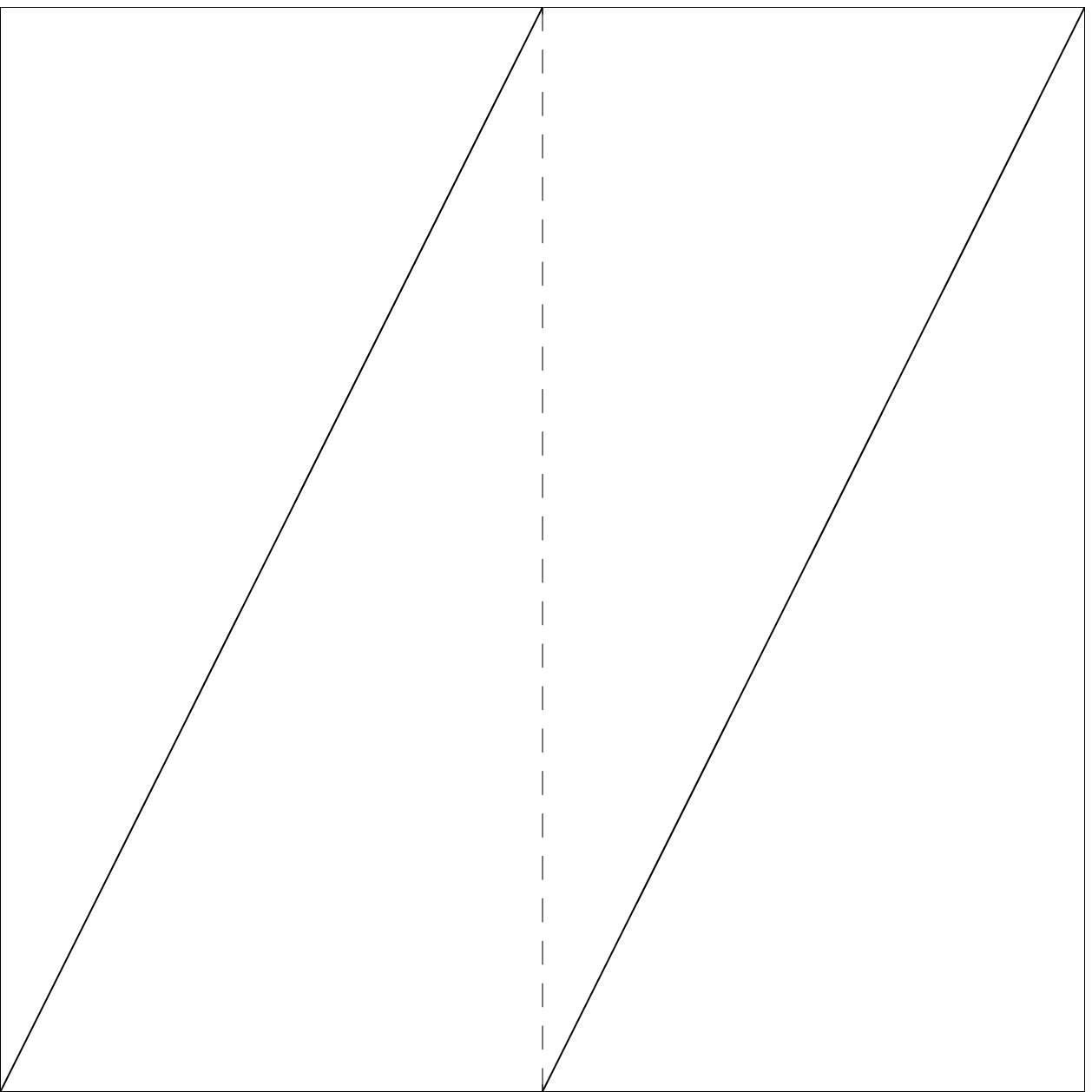}}
\put(-2,-2){\makebox(0,0)[tr]{$0$}}
\put(180,-2){\makebox(0,0)[t]{$\frac{1}{2}$}}
\put(360,-2){\makebox(0,0)[t]{$1$}}
\put(-2,360){\makebox(0,0)[r]{$1$}}
\end{picture}\kern60\unitlength\begin{picture}(390,390)(-15,-15)
\put(0,1){\includegraphics[bb=0 0 360 360,width=108bp,height=108bp]{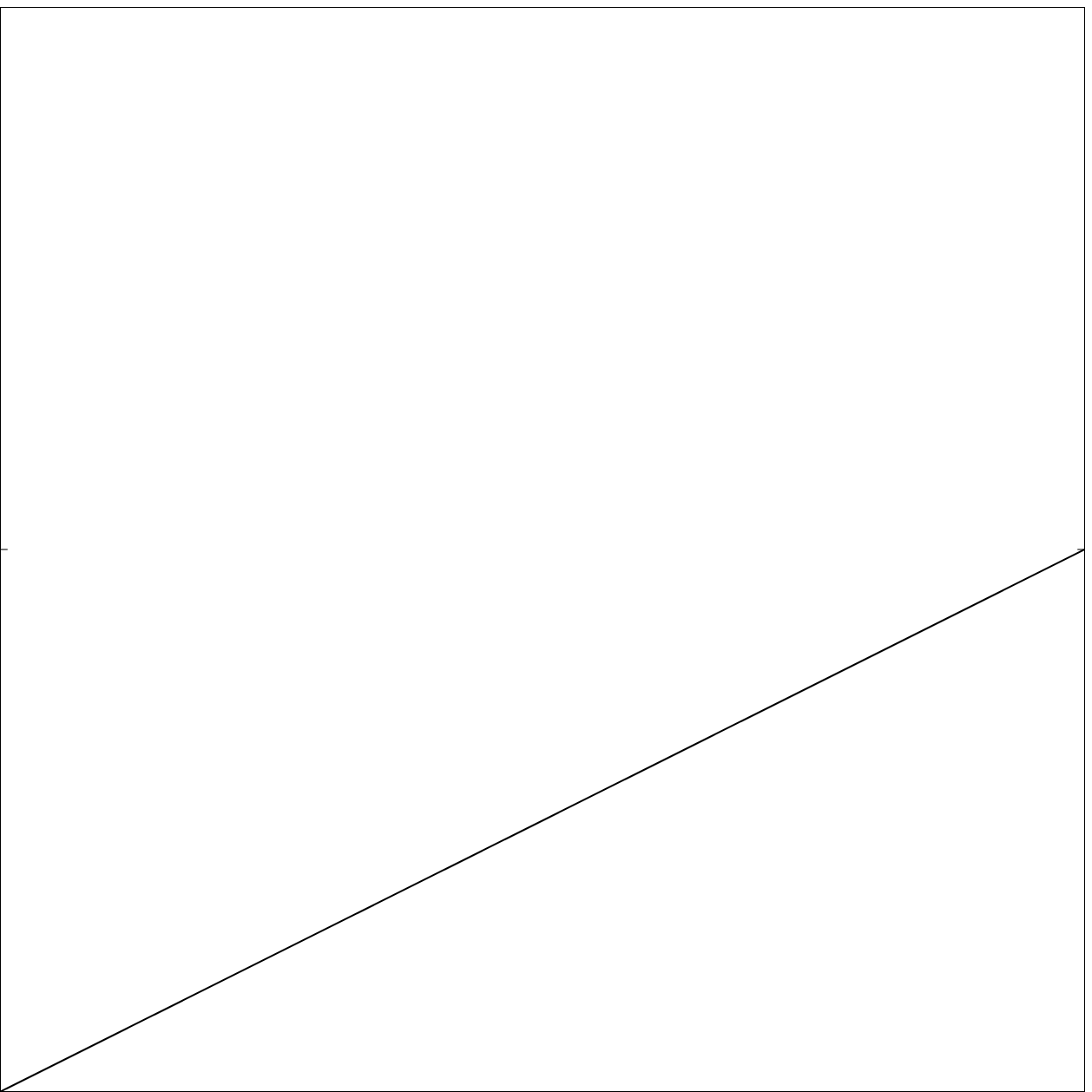}}
\put(-2,-2){\makebox(0,0)[tr]{$0$}}
\put(360,-2){\makebox(0,0)[t]{$1$}}
\put(-2,360){\makebox(0,0)[r]{$1$}}
\put(-2,180){\makebox(0,0)[r]{$\frac{1}{2}$}}
\end{picture}\kern60\unitlength\begin{picture}(390,390)(-15,-15)
\put(0,1){\includegraphics[bb=0 0 360 360,width=108bp,height=108bp]{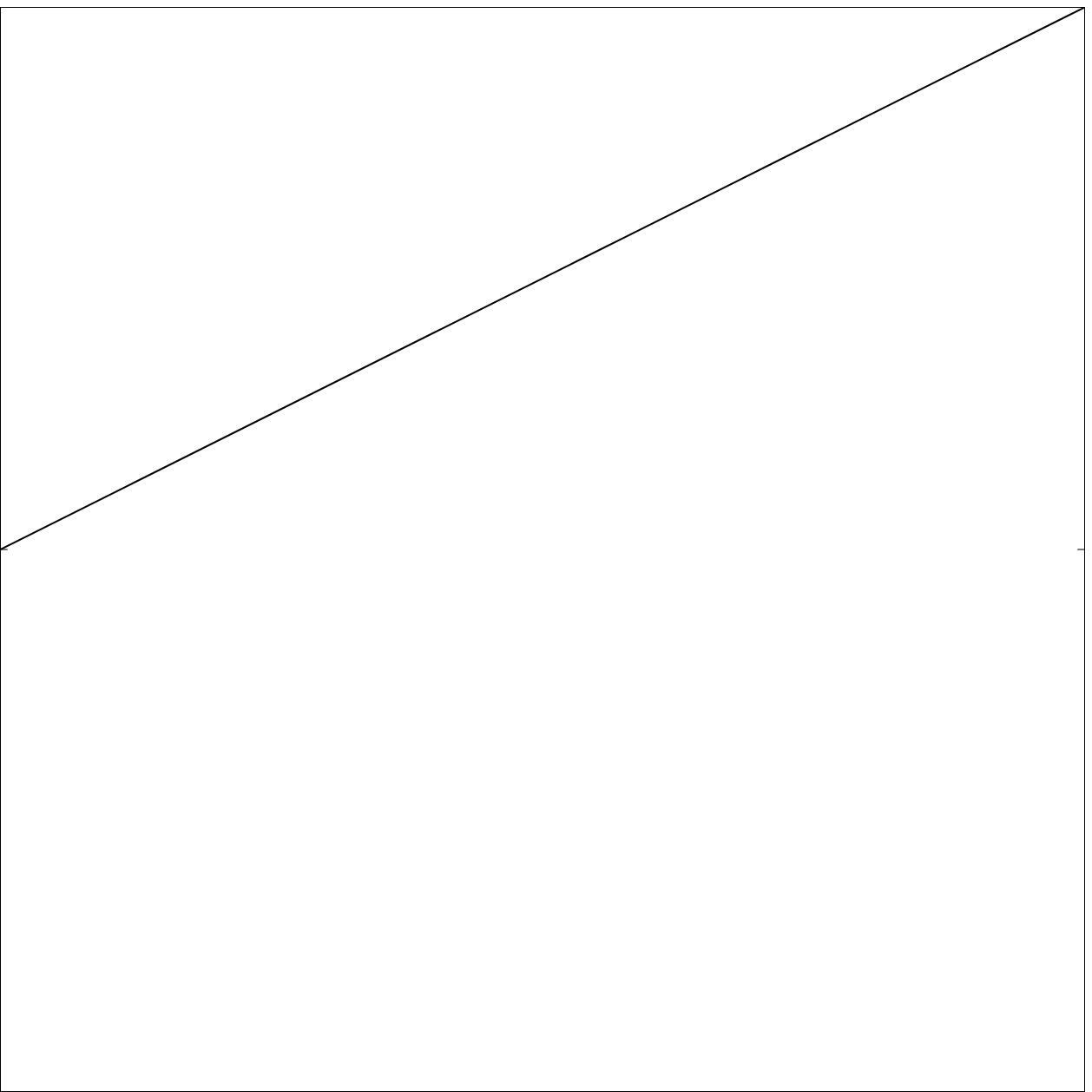}}
\put(-2,-2){\makebox(0,0)[tr]{$0$}}
\put(360,-2){\makebox(0,0)[t]{$1$}}
\put(-2,360){\makebox(0,0)[r]{$1$}}
\put(-2,180){\makebox(0,0)[r]{$\frac{1}{2}$}}
\end{picture}
\end{center}
\caption{Graphic representation of $\sigma$, $\sigma_{0}$, and $\sigma_{1}$.}%
\label{FigSigma01}%
\end{figure}

\begin{exmp}
\label{ExaSub.1}The particular example, $N=2$, mentioned above arises this way
when the identification
\begin{equation}
\left[  0,1\right)  \cong\mathbb{R}/\mathbb{Z} \label{eqSub.3}%
\end{equation}
is made and $\sigma\left(  x\right)  =2x\bmod{1}$. The three maps $\sigma$,
$\sigma_{0}$, and $\sigma_{1}$ may then be represented by the graphs in Figure
\ref{FigSigma01}.

If $J=\left[  0,1\right)  $ is the usual unit interval on the line, then the
subdivision from (\ref{eqSub.1}) takes the form%
\begin{equation}
J_{k}\left(  a\right)  =\sigma_{a_{1}}\circ\sigma_{a_{2}}\circ\dots\circ
\sigma_{a_{k}}\left(  J\right)  , \label{eqSub.4}%
\end{equation}
and the system%
\begin{equation}
\sigma_{a}:=\sigma_{a_{1}}\circ\sigma_{a_{2}}\circ\dots\circ\sigma_{a_{k}}
\label{eqSub.5}%
\end{equation}
forms a set of branches for $\sigma^{k}=\underbrace{\sigma\circ\dots
\circ\sigma}_{k\text{ times}}$ as $a_{i}\in\left\{  0,1\right\}  $,
$i=1,\dots,k$.
\end{exmp}

\begin{exmp}
\label{ExaSub.2}Set $\sigma_{0}\left(  x\right)  =x/3$ and $\sigma_{1}\left(
x\right)  =\left(  x+2\right)  /3$, and let $X\subset\mathbb{R}$ be the unique
solution to
\begin{equation}
X=\sigma_{0}\left(  X\right)  \cup\sigma_{1}\left(  X\right)  .
\label{eqSub.6}%
\end{equation}
Then $X$ is the familiar middle-third Cantor set, and there is a unique Borel
probability measure $\mu$ supported on $X$ and satisfying%
\begin{equation}
\mu=\frac{1}{2}\left(  \mu\circ\sigma_{0}^{-1}+\mu\circ\sigma_{1}^{-1}\right)
\label{eqSub.7}%
\end{equation}
or equivalently%
\begin{equation}
\int f\,d\mu=\frac{1}{2}\left(  \int f\circ\sigma_{0}\,d\mu+\int f\circ
\sigma_{1}\,d\mu\right)  \text{\qquad for all }f\in C\left(  X\right)  .
\label{eqSub.8}%
\end{equation}

\end{exmp}

\section{\label{Pro}Projection-valued measures}

In this section we study subdivisions of compact metric spaces, and
subdivisions of projections in Hilbert space, and we use our observations in
the construction of certain projection-valued measures.

\begin{defn}
\label{DefPro.I}Let $\left(  X,d\right)  $ be a compact metric space. For
subsets $A\subset X$, we define the \emph{diameter}%
\begin{equation}
\left\vert A\right\vert :=\sup\left\{  \,d\left(  x,y\right)  \mid x,y\in
A\,\right\}  . \label{eqPro.1}%
\end{equation}
A \emph{partition} of $X$ is a family $\left\{  A\left(  i\right)  \right\}
_{i\in I}$, $I$ some index set, such that%
\begin{equation}
\bigcup_{i}A\left(  i\right)  =X\text{,\quad and\quad}A\left(  i\right)  \cap
A\left(  j\right)  =\varnothing\text{ if }i\neq j. \label{eqPro.2}%
\end{equation}
Let $N\in\mathbb{Z}_{+}$, $N\geq2$. Let $\Gamma_{\!N}:=\left\{  0,1,\dots
,N-1\right\}  $. Suppose for each $k\in\mathbb{Z}_{+}$, we have a partition
into Borel subsets $\left\{  A_{k}\left(  a\right)  \right\}  $ indexed by
$a\in\Gamma_{\!N}^{k}=\underbrace{\Gamma_{\!N}\times\dots\times\Gamma_{\!N}%
}\limits_{k\text{ times}}$, and%
\begin{equation}
\left\vert A_{k}\left(  a\right)  \right\vert =O\left(  N^{-ck}\right)
,\qquad c>0. \label{eqPro.3}%
\end{equation}
Suppose every $A_{k+1}\left(  a\right)  $ is contained in some $A_{k}\left(
b\right)  $. We then say that $\left\{  A_{k}\left(  a\right)  \right\}  $ is
an $N$\emph{-adic system of partitions} of $X$.
\end{defn}

\begin{defn}
\label{DefPro.II}Let $\mathcal{H}$ be a complex Hilbert space. A
\emph{partition of projections} in $\mathcal{H}$ is a system $\left\{
P\left(  i\right)  \right\}  _{i\in I}$ of projections, i.e., $P\left(
i\right)  =P\left(  i\right)  ^{\ast}=P\left(  i\right)  ^{2}$, such that
\begin{equation}
P\left(  i\right)  P\left(  j\right)  =0\text{ if }i\neq j\text{\quad
and\quad}\sum_{i}P\left(  i\right)
=\hbox{\upshape \small1\kern-3.8pt\normalsize1}_{\mathcal{H}}, \label{eqPro.4}%
\end{equation}
where $\hbox{\upshape \small1\kern-3.8pt\normalsize1}_{\mathcal{H}}$ denotes
the identity operator in $\mathcal{H}$. Let $N\in\mathbb{Z}_{+}$, $N\geq2$.
Suppose for each $k\in\mathbb{Z}_{+}\ $we have a partition of projections
$\left\{  P_{k}\left(  a\right)  \right\}  _{a\in\Gamma_{\!N}^{k}}$ such that
every $P_{k+1}\left(  a\right)  $ is contained in some $P_{k}\left(  b\right)
$, i.e.,%
\begin{equation}
P_{k}\left(  b\right)  P_{k+1}\left(  a\right)  =P_{k+1}\left(  a\right)  .
\label{eqPro.5}%
\end{equation}
Then the combined system $\left\{  P_{k}\left(  a\right)  \right\}
_{k\in\mathbb{Z}_{+},~a\in\Gamma_{\!N}^{k}}$ is a system of partitions of
$\hbox{\upshape \small1\kern-3.8pt\normalsize1}_{\mathcal{H}}$ forming, by
$N$-adic subdivisions, an $N$-adic system of projections. We refer to this
system as an $N$\emph{-adic system of partitions of }%
$\hbox{\upshape \small1\kern-3.8pt\normalsize1}_{\mathcal{H}}$ \emph{into
projections} (or for brevity, if the context makes clear what Hilbert space
$\mathcal{H}$ is being partitioned into projections, simply as an
$N$\emph{-adic system of projections\/}).
\end{defn}

\begin{defn}
\label{RemPro.1}We use $\mathcal{B}\left(  X\right)  $ to denote the Borel
subsets of the compact metric space $X$. A positive operator-valued function
$E$ defined on $\mathcal{B}\left(  X\right)  $ may be called a $\sigma
$\emph{-additive measure} if, given a sequence $B_{1},B_{2},\dots$ in
$\mathcal{B}\left(  X\right)  $ such that $B_{i}\cap B_{j}=\varnothing$ for
$i\neq j$, the measures combine according to the formula%
\begin{equation}
E\left(  
\bigcup_{i}B_{i}
\right)  =\sum_{i}E\left(  B_{i}\right)  ,
\label{eqPro.7}%
\end{equation}
where, since the values $E\left(  B_{i}\right)  $ are positive operators, we
may take the summation on the right-hand side to be convergent in the strong
operator topology. Such a measure is called an \emph{orthogonal
projection-valued measure} if it satisfies the additional properties
(\ref{LemPro.1(1)})--(\ref{LemPro.1(3)}):

\begin{enumerate}
\item \label{LemPro.1(1)}$E\left(  B\right)  =E\left(  B\right)  ^{\ast
}=E\left(  B\right)  ^{2}$ for $B\in\mathcal{B}\left(  X\right)  ={}$the Borel
subsets of $X$,

\item \label{LemPro.1(2)}$E\left(  B_{1}\right)  E\left(  B_{2}\right)  =0$ if
$B_{1},B_{2}\in\mathcal{B}\left(  X\right)  $ satisfy $B_{1}\cap
B_{2}=\varnothing$, and

\item \label{LemPro.1(3)}$E\left(  X\right)
=\hbox{\upshape \small1\kern-3.8pt\normalsize1}_{\mathcal{H}}$.
\end{enumerate}
\end{defn}

\begin{rem}
\label{RemPro.pound}There are four independent conditions in Definition
\ref{RemPro.1}. If $E\left(  \,\cdot\,\right)  $ is a function defined on
$\mathcal{B}\left(  X\right)  $ and taking values in positive operators on
$\mathcal{H}$, and if only property (\ref{eqPro.7}) is satisfied, we say
$E\left(  \,\cdot\,\right)  $ is a positive operator-valued measure. If
(\ref{LemPro.1(1)}) is \emph{also} satisfied, we say that $E$ is
projection-valued. If (\ref{eqPro.7}), (\ref{LemPro.1(1)}), and
(\ref{LemPro.1(2)}) are satisfied, we say that the projection-valued measure
is orthogonal. If all four conditions hold, we talk of a projection-valued
measure which is orthogonal and normalized, in short an \emph{orthogonal
projection-valued measure}. In this paper, we will only have occasion to study
the case when $E\left(  \,\cdot\,\right)  $ satisfies all four conditions.
\end{rem}

\begin{lem}
\label{LemPro.1}Let $N\in\mathbb{Z}_{+}$, $N\geq2$. Let $\left(  X,d\right)  $
be a compact metric space, and let $\mathcal{H}$ be a complex Hilbert space.
Let $\left\{  A_{k}\left(  a\right)  \right\}  _{k\in\mathbb{Z}_{+}%
,\,a\in\Gamma_{\!N}^{k}}$ be an $N$-adic system of partitions of $X$, and let
$\left\{  P_{k}\left(  a\right)  \right\}  _{k\in\mathbb{Z}_{+},\,a\in
\Gamma_{\!N}^{k}}$ be the corresponding $N$-adic system of partitions of
$\hbox{\upshape \small1\kern-3.8pt\normalsize1}_{\mathcal{H}}$ into
projections \textup{(}i.e., the corresponding $N$-adic system of
projections\/\textup{)}. Then there is a unique normalized orthogonal
projection-valued measure $E\left(  \,\cdot\,\right)  $ defined on the Borel
subsets of $X$ and taking values in the orthogonal projections in
$\mathcal{H}$ such that%
\begin{equation}
E\left(  A_{k}\left(  a\right)  \right)  =P_{k}\left(  a\right)  \text{\qquad
for all }k\in\mathbb{Z}_{+}\text{ and }a\in\Gamma_{\!N}^{k}. \label{eqPro.6}%
\end{equation}

\end{lem}

\begin{defn}
\label{DefPro.III}Let $N\in\mathbb{Z}_{+}$, $N\geq2$. We shall need the
\emph{Cuntz algebra} $\mathcal{O}_{N}$ \cite{Cun77} on $N$ generators
$s_{0},s_{1},\dots,s_{N-1}$. It is the unique $C^{\ast}$-algebra on the
relations%
\begin{equation}
\sum_{i}s_{i}s_{i}^{\ast}=\hbox{\upshape \small1\kern-3.8pt\normalsize1}
\label{eqPro.8}%
\end{equation}
where $\hbox{\upshape \small1\kern-3.8pt\normalsize1}$ is the unit element in
the algebra $\mathcal{O}_{N}$. We have
\begin{equation}
s_{i}^{\ast}s_{j}=\delta_{i,j}\hbox{\upshape \small1\kern-3.8pt\normalsize1}.
\label{eqPro.9}%
\end{equation}
To specify a \emph{representation of} $\mathcal{O}_{N}$ on a Hilbert space
$\mathcal{H}$ we need $N$ isometries $S_{0},S_{1},\dots,S_{N-1}$ such that%
\begin{equation}
\sum_{i}S_{i}S_{i}^{\ast}%
=\hbox{\upshape \small1\kern-3.8pt\normalsize1}_{\mathcal{H}}.
\label{eqPro.10}%
\end{equation}
Then $S_{i}^{\ast}S_{j}=\delta_{i,j}%
\hbox{\upshape \small1\kern-3.8pt\normalsize1}_{\mathcal{H}}$ and the
representation is determined uniquely.
\end{defn}

\begin{lem}
\label{LemPro.2}Let $N\in\mathbb{Z}_{+}$, $N\geq2$, and let $S_{0},S_{1}%
,\dots,S_{N-1}$ be a representation of $\mathcal{O}_{N}$ on a Hilbert space
$\mathcal{H}$. For $k\in\mathbb{Z}_{+}$ and $a=\left(  a_{1},\dots
,a_{k}\right)  \in\Gamma_{\!N}^{k}$, set
\begin{equation}
S_{a}:=S_{a_{1}}\cdots S_{a_{k}}\text{\quad and\quad}P_{k}\left(  a\right)
=S_{a}S_{a}^{\ast}. \label{eqPro.11}%
\end{equation}
Then the combined system $\left\{  P_{k}\left(  a\right)  \right\}
_{k\in\mathbb{Z}_{+},\,a\in\Gamma_{\!N}^{k}}$ is an $N$-adic system of
partitions of $\hbox{\upshape \small1\kern-3.8pt\normalsize1}_{\mathcal{H}}$
into projections \textup{(}i.e., an $N$-adic system of projections\/\textup{)}.
\end{lem}

We now turn to the proof of the two lemmas.%

\begin{pf*}{PROOF OF LEMMA \ref{LemPro.1}.}%
Let $N\in\mathbb{Z}_{+}$, and let systems $\left\{  A_{k}\left(  a\right)
\right\}  _{k\in\mathbb{Z}_{+},\,a\in\Gamma_{\!N}^{k}}$ and $\left\{
P_{k}\left(  a\right)  \right\}  _{k\in\mathbb{Z}_{+},\,a\in\Gamma_{\!N}^{k}}$
be given as in the statement of the lemma. For every $k\in\mathbb{Z}_{+}$, the
finite sums%
\begin{equation}
\sum_{a\in\Gamma_{\!N}^{k}}C_{a}\chi_{A_{k}\left(  a\right)  }\label{eqPro.12}%
\end{equation}
form an algebra $\mathfrak{A}_{k}$ of functions on $X$, and from the
definition of the partition system $\left\{  A_{k}\left(  a\right)  \right\}
_{k\in\mathbb{Z}_{+},\,a\in\Gamma_{\!N}^{k}}$ it follows that there are
natural embeddings $\mathfrak{A}_{k}\hookrightarrow\mathfrak{A}_{k+1}$. {}From
the definition of the projection system $\left\{  P_{k}\left(  a\right)
\right\}  _{k\in\mathbb{Z}_{+},\,a\in\Gamma_{\!N}^{k}}$ we conclude that the
mappings, defined for each $k\in\mathbb{Z}_{+}$,%
\begin{equation}
\sum_{a\in\Gamma_{\!N}^{k}}C_{a}\chi_{A_{k}\left(  a\right)  }\longmapsto
\sum_{a\in\Gamma_{\!N}^{k}}C_{a}P_{k}\left(  a\right)  \label{eqPro.13}%
\end{equation}
extends to the algebra%
\begin{equation}
\mathfrak{A}:=\bigcup_{\phantom{{}_{+}}k\in\mathbb{Z}_{+}}\mathfrak{A}%
_{k}.\label{eqPro.14}%
\end{equation}
But the operators on the right-hand side in (\ref{eqPro.13}) form an abelian
algebra $\mathcal{C}$ of operators. The algebra $\mathfrak{A}$ is closed under
complex conjugation $f\mapsto\bar{f}$, and $\mathcal{C}$ is $\ast$-closed,
i.e., $E\in\mathcal{C}\Rightarrow E^{\ast}\in\mathcal{C}$. Let the mapping
obtained from (\ref{eqPro.13}) be denoted $\pi$. Then one checks from the two
definitions, Definition \ref{DefPro.I} and Definition \ref{DefPro.II}, that
$\pi\left(  f_{1}f_{2}\right)  =\pi\left(  f_{1}\right)  \pi\left(
f_{2}\right)  $, $f_{1},f_{2}\in\mathfrak{A}$, and $\pi\left(
\smash{\bar{f}}\vphantom{f}\,\right)  =\pi\left(  f\right)  ^{\ast}$,
$f\in\mathfrak{A}$. Since the sets $A_{k}\left(  a\right)  $ satisfy
$\left\vert A_{k}\left(  a\right)  \right\vert =O\left(  N^{-k}\right)  $,
where $\left\vert \,\cdot\,\right\vert $ denotes the diameter, it is clear
that every $f\in C\left(  X\right)  $ may be approximated uniformly with a
sequence in $\mathfrak{A}$. This means that $f\mapsto\pi\left(  f\right)  $
extends uniquely $\tilde{\pi}\colon C\left(  X\right)  \rightarrow B\left(
\mathcal{H}\right)  $, and the extension satisfies%
\begin{equation}
\tilde{\pi}\left(  f_{1}f_{2}\right)  =\tilde{\pi}\left(  f_{1}\right)
\tilde{\pi}\left(  f_{2}\right)  ,\;f_{1},f_{2}\in C\left(  X\right)
\text{,\quad and\quad}\tilde{\pi}\left(  \smash{\bar{f}}\vphantom{f}\,\right)
=\tilde{\pi}\left(  f\right)  ^{\ast}.\label{eqPro.15}%
\end{equation}
A standard argument from function theory now shows that $\tilde{\pi}$ extends
further from $C\left(  X\right)  $ to all the Baire functions, and the
extension satisfies the same multiplication rules (\ref{eqPro.15}). For this
part of the argument see, e.g., \cite[Section 6]{Nel69}. If $B\in
\mathcal{B}\left(  X\right)  $, we may then define a projection-valued measure
$E\left(  \,\cdot\,\right)  $ as follows:%
\begin{equation}
E\left(  B\right)  :=\tilde{\pi}\left(  \chi_{B}\right)  ,\label{eqPro.16}%
\end{equation}
where $\chi_{B}$ denotes the indicator function of the set $B$. Since
$\tilde{\pi}$ is obtained as a unique extension from (\ref{eqPro.13}) it
follows immediately that $E\left(  \,\cdot\,\right)  $ in (\ref{eqPro.16}) has
the properties (\ref{LemPro.1(1)})--(\ref{LemPro.1(3)}) from Definition
\ref{RemPro.1}, and that it is countably additive, see (\ref{eqPro.7}).
Moreover, it satisfies (\ref{eqPro.6}), and is determined uniquely by
(\ref{eqPro.6}). \qed%
\end{pf*}%
%

\begin{pf*}{PROOF OF LEMMA \ref{LemPro.2}.}%
The details are essentially well known; see, e.g., \cite{BJO04}. In fact, an
inspection shows that the projections $P_{k}\left(  a\right)  =S_{a}%
S_{a}^{\ast}$, $a\in\Gamma_{\!N}^{k}$, introduced in (\ref{eqPro.11}) generate
an abelian algebra of operators. It is a special case of the algebra
$\mathcal{C}$ which we introduced in the proof of Lemma \ref{LemPro.1} above.
Also note the following monotonicity: If $S$ and $T$ are positive operators on
$\mathcal{H}$, we say that $S\leq T$ if
\begin{equation}
\left\langle \,x\mid Sx\,\right\rangle \leq\left\langle \,x\mid
Tx\,\right\rangle \text{\quad holds for all }x\in\mathcal{H}.\label{eqPro.17}%
\end{equation}
The inner product of $\mathcal{H}$ is denoted $\left\langle \,\cdot\mid
\cdot\,\right\rangle $ and is assumed linear in the second factor. Using the
defining relation (\ref{eqPro.10}) for the generators of a representation of
$\mathcal{O}_{N}$, note that if $a\in\Gamma_{\!N}^{k}$ for some $k$, and if
$i\in\Gamma_{\!N}$, then
\begin{equation}
\left(  ai\right)  \in\Gamma_{\!N}^{k+1}\text{\quad and\quad}\sum_{i}%
P_{k+1}\left(  ai\right)  =P_{k}\left(  a\right)  .\label{eqPro.18}%
\end{equation}
As a result, we get $P_{k+1}\left(  ai\right)  \leq P_{k}\left(  a\right)  $,
or equivalently $P_{k}\left(  a\right)  P_{k+1}\left(  ai\right)
=P_{k+1}\left(  ai\right)  $, which is the desired relation (\ref{eqPro.5})
from Definition \ref{DefPro.II}. \qed%
\end{pf*}%

\begin{rem}
\label{RemPro.2}If $\left\{  S_{i}\right\}  _{i\in\Gamma_{\!N}}$ is a
representation of $\mathcal{O}_{N}$ on a Hilbert space, it is known that the
$C^{\ast}$-algebra $\mathcal{C}$ generated by the projections $P_{k}\left(
a\right)  =S_{a}S_{a}^{\ast}$, $a\in\Gamma_{\!N}^{k}$, $k\in\mathbb{Z}_{+}$,
is naturally isomorphic to the algebra of all continuous functions on the
infinite Cartesian
\begin{equation}
\text{product\quad}\prod_{\phantom{{}_{+}}\mathbb{Z}_{+}}\Gamma_{\!N}%
\text{\quad or\quad}\Gamma_{\!N}^{\mathbb{Z}_{+}} \label{eqPro.19}%
\end{equation}
when the Cartesian product is equipped with the product topology of Tychonoff,
i.e., $\mathcal{C}\cong C\left(  X\right)  $ with $X=\Gamma_{\!N}%
^{\mathbb{Z}_{+}}$. Recall that $X$ is compact; see, e.g., \cite{Rud91}.

It is easy to see that $X=\Gamma_{\!N}^{\mathbb{Z}_{+}}$ carries a system of
functions $\sigma,\sigma_{0},\dots,\sigma_{N-1}$ with the properties listed in
(\ref{eqSub.2a}). If elements $x$ in $X$ are represented as sequences $\left(
x_{1},x_{2},\dots\right)  $, we set
\begin{equation}
\sigma\left(  x_{1},x_{2},\dots\right)  =\left(  x_{2},x_{3},\dots\right)
\text{\quad and\quad}\sigma_{i}\left(  x_{1},x_{2},\dots\right)  =\left(
i,x_{1},x_{2},\dots\right)  . \label{eqPro.20}%
\end{equation}
All $N+1$ functions are continuous $X\rightarrow X$ and satisfy
(\ref{eqSub.2a}), i.e.,
\begin{equation}
\sigma\left(  \sigma_{i}\left(  x\right)  \right)  =x\text{\qquad for all
}x\in X\text{ and }i=0,1,\dots,N-1. \label{eqPro.21}%
\end{equation}

There is a family of measures on $X$ which generalizes the property
(\ref{eqSub.7}) above. They are the product measures: Let $\left\{
p_{i}\right\}  _{i\in\Gamma_{\!N}}$ be given such that $p_{i}\geq0$ and
$\sum_{i}p_{i}=1$, and let $\xi_{k}\colon X\rightarrow\Gamma_{\!N}$ be the
coordinate projection $\xi_{k}\left(  x_{1},x_{2},\dots\right)  =x_{k}$. For
subsets $T\subset\Gamma_{\!N}$, set%
\begin{equation}
\xi_{k}^{-1}\left(  T\right)  =\left\{  \, x\in X\mid\xi_{k}\left(  x\right)
\in T\,\right\}  . \label{eqPro.22}%
\end{equation}
Using standard measure theory \cite{Rud91}, note that there is a unique
measure $\mu_{p}$ on $X$ such that%
\begin{equation}
\mu_{p}\left(  \xi_{k}^{-1}\left(  \left\{  i\right\}  \right)  \right)
=p_{i}\text{\qquad for all }k\in\mathbb{Z}_{+},\;i\in\Gamma_{\!N}.
\label{eqPro.23}%
\end{equation}
Introducing the maps $\sigma_{i}\colon X\rightarrow X$ of (\ref{eqPro.20}) we
note that $\mu_{p}$ satisfies%
\begin{equation}
\mu_{p}=\sum_{i\in\Gamma_{\!N}}p_{i}\mu_{p}\circ\sigma_{i}^{-1},
\label{eqPro.24}%
\end{equation}
or equivalently%
\begin{equation}
\int_{X}f\,d\mu_{p}=\sum_{i\in\Gamma_{\!N}}p_{i}\int_{X}f\circ\sigma_{i}%
\,d\mu_{p}\text{\qquad for all }f\in C\left(  X\right)  . \label{eqPro.25}%
\end{equation}
Finally, note that distinct probabilities $\left(  p_{i}\right)  $ and
$\left(  p_{i}^{\prime}\right)  $ yield measures $\mu_{p}$ and $\mu
_{p^{\prime}}$ which are mutually singular.

If $\left(  S_{i}\right)  _{0\leq i<N}$ is a representation of $\mathcal{O}%
_{N}$ for some $N\in\mathbb{Z}_{+}$, $N\geq2$, then we will denote the
corresponding projection-valued measure on $\Gamma_{\!N}^{\mathbb{Z}_{+}}$ by
$E$. If $\left\{  A_{k}\left(  a\right)  \right\}  _{k\in\mathbb{Z}_{+}%
,\,a\in\Gamma_{\!N}^{k}}$ is an $N$-adic system of partitions of some compact
metric space $\left(  X,d\right)  $, then the corresponding projection-valued
measure on $\mathcal{B}\left(  X\right)  $ will be denoted $E^{A}\left(
\,\cdot\,\right)  $ to stress its dependence on the partition system.
\end{rem}

The next lemma shows that the algebra $\mathcal{C}$ in the proof of Lemma
\ref{LemPro.1} is isomorphic to $C\left(  \Gamma_{\!N}^{\mathbb{Z}_{+}%
}\right)  $ where $\mathcal{C}$ is viewed as a $C^{\ast}$-algebra, and the
infinite Cartesian product $X=\Gamma_{\!N}^{\mathbb{Z}_{+}}$ is given its
Tychonoff topology. Since $\mathcal{C}$ is an abelian $C^{\ast}$-algebra we
know that it is isomorphic to $C\left(  K\right)  $ for some compact Hausdorff
space $K$. The isomorphism $\mathcal{C}\cong C\left(  K\right)  $ is called
the Gelfand transform, and $K$ the Gelfand space. The conclusion of the lemma
is that $\Gamma_{\!N}^{\mathbb{Z}_{+}}$ is the Gelfand space of $\mathcal{C}$,
and further we offer a formula for the Gelfand transform. We also note, by
standard theory, see, e.g., \cite{Rud91}, that $\Gamma_{2}^{\mathbb{Z}_{+}}$
is homeomorphic to the Cantor set $X$ in Example \ref{ExaSub.2} above. In
particular, it is known that the compact space $\Gamma_{\!N}^{\mathbb{Z}_{+}}$
is totally disconnected.

\begin{lem}
\label{LemPro.3}Let $N\in\mathbb{Z}_{+}$, $N\geq2$, and let $\mathcal{O}_{N}$
be the Cuntz $C^{\ast}$-algebra with generators $\left\{  s_{i}\right\}
_{0\leq i<N}$ subject to the axioms in Definition \textup{\ref{DefPro.III}},
i.e., \textup{(\ref{eqPro.8})}. The $C^{\ast}$-algebra $\mathcal{C}$ is the
norm-closure of the algebra generated by the elements%
\begin{equation}
e\left(  a\right)  :=s_{a}s_{a}^{\ast}=s_{a_{1}}s_{a_{2}}\cdots s_{a_{k}%
}s_{a_{k}}^{\ast}\cdots s_{a_{2}}^{\ast}s_{a_{1}}^{\ast}, \label{eqPro.26}%
\end{equation}
where $k\in\mathbb{Z}_{+}$ and $a=\left(  a_{1},\dots,a_{k}\right)  \in
\Gamma_{\!N}^{k}$. Let $\xi_{i}\colon X\rightarrow\Gamma_{\!N}$ be the
coordinate function \textup{(\ref{eqPro.22})}, $X=\Gamma_{\!N}^{\mathbb{Z}%
_{+}}$. Then the assignment
\begin{equation}
\left(  \chi_{\left\{  a_{1}\right\}  }\circ\xi_{1}\right)  \left(
\chi_{\left\{  a_{2}\right\}  }\circ\xi_{2}\right)  \cdots\left(
\chi_{\left\{  a_{k}\right\}  }\circ\xi_{k}\right)  \overset{G}{\longmapsto
}e\left(  a_{1},\dots,a_{k}\right)  \label{eqPro.27}%
\end{equation}
extends to a $C^{\ast}$-isomorphism of $C\left(  X\right)  $ onto
$\mathcal{C}$.
\end{lem}

\begin{pf}
The function $f_{a}$, for $a\in\Gamma_{\!N}^{k}$, in the formula on the
left-hand side in (\ref{eqPro.27}) is given as follows: Evaluation at
$x=\left(  x_{i}\right)  \in X$, $f_{a}\left(  x\right)  =\delta_{a_{1},x_{1}%
}\delta_{a_{2},x_{2}}\cdots\delta_{a_{k},x_{k}}$. {}From the definition of the
Tychonoff topology it follows that each $f_{a}$ is continuous, and that the
family $\left\{  f_{a}\right\}  _{k\in\mathbb{Z}_{+},\,a\in\Gamma_{\!N}^{k}}$
separates points in $X$.

It follows from the relations on the generators $\left\{  s_{i}\right\}
_{0\leq i<N}$ that the assignment $G$ in (\ref{eqPro.27}) is an isomorphism
from an abelian subalgebra $\mathcal{S}$ of $C\left(  X\right)  $ into a dense
subalgebra of $\mathcal{C}$. But $\mathcal{S}$ is dense in $C\left(  X\right)
$ by virtue of the Stone-Weierstra\ss \ theorem, and it is immediate from this
that $G$ extends uniquely, by closure, to a $\ast$-isomorphism of $C\left(
X\right)  $ onto $\mathcal{C}$. \qed
\end{pf}

\begin{defn}
\label{DefPro.4}Let $\left(  X,d\right)  $ be a compact metric space, and let
$N\in\mathbb{Z}_{+}$, $N\geq2$, be given. We say that an $N$-adic system
$\left\{  A_{k}\left(  a\right)  \right\}  _{k\in\mathbb{Z}_{+},\,a\in
\Gamma_{\!N}^{k}}$ of partitions of $X$ is \emph{affiliated with} an iterated
function system (IFS) on $X$ if there is a system $\sigma,\left(  \sigma
_{i}\right)  _{0\leq i<N}$ of continuous maps such that%
\begin{equation}
\sigma\circ\sigma_{i}=\operatorname*{id}\nolimits_{X},\qquad i\in\Gamma_{\!N},
\label{eqProNew.1}%
\end{equation}
and%
\begin{equation}
\sigma_{a_{1}}\sigma_{a_{2}}\cdots\sigma_{a_{k}}\left(  X\right)
=A_{k}\left(  a\right)  \text{\qquad for all }k\in\mathbb{Z}_{+}\text{ and
}a\in\Gamma_{\!N}^{k}. \label{eqProNew.2}%
\end{equation}
(Note that (\ref{eqProNew.1}) is part of the definition of an IFS.)
\end{defn}

The following is a corollary to the result stated as Lemma \ref{LemPro.1}
above, i.e., the construction of a projection-valued measure $E^{A}\left(
\,\cdot\,\right)  $ from a given representation $\left(  S_{i}\right)  $ of
$\mathcal{O}_{N}$ on a Hilbert space and a given $N$-adic system $\left(
A_{k}\left(  a\right)  \right)  $ of partitions.

\begin{cor}
\label{CorProNew.1}Let $\left(  S_{i}\right)  $ be a representation of
$\mathcal{O}_{N}$ on a Hilbert space $\mathcal{H}$, and let $\left(
A_{k}\left(  a\right)  \right)  $ be an $N$-adic system of partitions which is
affiliated with a continuous iterated function system $\sigma,\sigma_{i}$
acting on a compact metric space $\left(  X,d\right)  $. Then the
projection-valued measure $E^{A}\left(  \,\cdot\,\right)  $, see Lemma
\textup{\ref{LemPro.1}}, satisfies
\begin{multline}
\qquad S_{a}E^{A}\left(  B\right)  S_{a}^{\ast}=E^{A}\left(  \sigma
_{a}B\right)  \text{\qquad for all }B\in\mathcal{B}\left(  X\right)
\text{,}\label{eqProNew.3}\\
k\in\mathbb{Z}_{+}\text{, and all }a\in\Gamma_{\!N}^{k},\qquad
\end{multline}
and%
\begin{equation}
\sum_{i=0}^{N-1}S_{i}E^{A}\left(  B\right)  S_{i}^{\ast}=E^{A}\left(
\sigma^{-1}\left(  B\right)  \right)  , \label{eqProNew.4}%
\end{equation}
where $a=\left(  a_{1},\dots,a_{k}\right)  $, $\sigma_{a}=\sigma_{a_{1}}%
\circ\dots\circ\sigma_{a_{k}}$, and%
\[
\sigma^{-1}\left(  B\right)  =\left\{  \, x\in X\mid\sigma\left(  x\right)
\in B\,\right\}  .
\]

\end{cor}

\begin{pf}
The argument in the proof of (\ref{eqProNew.3}) and (\ref{eqProNew.4}) is
based directly on the two-step approximation which went into the construction
of the measure $E^{A}\left(  \,\cdot\,\right)  $; see Lemma \ref{LemPro.1} for details.

With the assumptions on the representation $\left(  S_{i}\right)  $ and the
IFS partition, the two operator commutation relations (\ref{eqProNew.3}%
)--(\ref{eqProNew.4}) follow from the same approximation, coupled with the
observation that if $a\in\Gamma_{\!N}^{k}$ and $b\in\Gamma_{\!N}^{l}$, then
\begin{align}
S_{a}E^{A}\left(  A_{l}\left(  b\right)  \right)  S_{a}^{\ast} &  =S_{a}%
S_{b}S_{b}^{\ast}S_{a}^{\ast}=S_{ab}S_{ab}^{\ast}\nonumber\\
&  =E^{A}\left(  A_{k+l}\left(  ab\right)  \right)  =E^{A}\left(  \sigma
_{a}\left(  A_{l}\left(  b\right)  \right)  \right)  ,\label{eqProNew.5}%
\end{align}
where $ab=\left(  a_{1},\dots,a_{k},b_{1},\dots,b_{l}
\right)  \in\Gamma_{\!N}^{k+1}$, i.e., concatenation, and the formula
\begin{equation}
\bigcup_{i}\sigma_{i}B=\sigma^{-1}\left(  B\right)  ,\qquad B\in
\mathcal{B}\left(  X\right)  .\label{eqProFinal}%
\end{equation}
\qed
\end{pf}

\section{\label{End}Endomorphisms of $B\left(  \mathcal{H}\right)  $}

{}From the point of view of the pure mathematics of operator algebras, it is
natural to ask about the von Neumann type
of the representations of the Cuntz algebras
that come from subband filters
(i.e., are they type 
$\mathrm{I}$, $\mathrm{II}$, or $\mathrm{III}$,
and how do they decompose?)
While this is addressed in \cite{Jor01b},
and to some extent (in a different context) in 
\cite{BJO04}, \cite{BJP96}, and \cite{BrJo99a}, we will not
discuss it here. Rather we will address a related question regarding the
selection of the ``best'' wavelets in specific parametrized families.

Let $N\in\mathbb{Z}_{+}$, $N\geq2$, and let $\mathcal{H}$ be a complex Hilbert
space. In our understanding of scaling problems in approximation theory, it is
often helpful to study endomorphisms of the $C^{\ast}$-algebra of all bounded
operators on $\mathcal{H}$. By this we mean a linear mapping $\alpha\colon
B\left(  \mathcal{H}\right)  \rightarrow B\left(  \mathcal{H}\right)  $ taking
$\hbox{\upshape \small1\kern-3.8pt\normalsize1}_{\mathcal{H}}$ to
$\hbox{\upshape \small1\kern-3.8pt\normalsize1}_{\mathcal{H}}$ and satisfying
\begin{equation}
\alpha\left(  ST\right)  =\alpha\left(  S\right)  \alpha\left(  T\right)
,\;S,T\in B\left(  \mathcal{H}\right)  ,\text{\quad and\quad}\alpha\left(
T^{\ast}\right)  =\alpha\left(  T\right)  ^{\ast}. \label{eqEnd.1}%
\end{equation}
We showed in \cite{BJP96} that there is a correspondence between
$\operatorname*{End}\left(  B\left(  \mathcal{H}\right)  \right)  $ and
representations of the Cuntz relations; see (\ref{eqPro.8}) in Definition
\ref{DefPro.III}. If $\left(  S_{i}\right)  _{0\leq i<N}$ is a representation
of $\mathcal{O}_{N}$ on a Hilbert space $\mathcal{H}$, then define
$\alpha\colon B\left(  \mathcal{H}\right)  \rightarrow B\left(  \mathcal{H}%
\right)  $ by%
\begin{equation}
\alpha\left(  T\right)  =\sum_{i=0}^{N-1}S_{i}TS_{i}^{\ast},\qquad T\in
B\left(  \mathcal{H}\right)  , \label{eqEnd.2}%
\end{equation}
and it is clear that $\alpha\in\operatorname*{End}\left(  B\left(
\mathcal{H}\right)  \right)  $. Moreover, the relative commutant%
\[
B\left(  \mathcal{H}\right)  \cap\alpha\left(  B\left(  \mathcal{H}\right)
\right)  ^{\prime}=\left\{  \, T\in B\left(  \mathcal{H}\right)  \mid
T\alpha\left(  X\right)  =\alpha\left(  X\right)  T,\;X\in B\left(
\mathcal{H}\right)  \,\right\}
\]
is naturally isomorphic to the algebra of all $N$-by-$N$ complex matrices
$M_{N}\left(  \mathbb{C}\right)  $. If $\left(  e_{i,j}\right)  $ are the
usual matrix units in $M_{N}\left(  \mathbb{C}\right)  $, i.e.,%
\begin{equation}
e_{i,j}\left(  k,l\right)  =\delta_{i,k}\delta_{j,l}\,, \label{eqEnd.3}%
\end{equation}
then the assignment
\begin{equation}
M_{N}\left(  \mathbb{C}\right)  \ni e_{i,j}\longmapsto S_{i}S_{j}^{\ast}%
\in\alpha\left(  B\left(  \mathcal{H}\right)  \right)  ^{\prime}
\label{eqEnd.4}%
\end{equation}
defines the isomorphism. There is a similar result for the commutant of the
iterated mapping $\alpha^{k}\left(  B\left(  \mathcal{H}\right)  \right)  $.
Then there is a similar isomorphism:%
\begin{equation}
e_{i_{1},j_{1}}\otimes e_{i_{2},j_{2}}\otimes\dots\otimes e_{i_{k},j_{k}%
}\longmapsto S_{i_{1}}S_{i_{2}}\cdots S_{i_{k}}S_{j_{k}}^{\ast}\cdots
S_{j_{2}}^{\ast}S_{j_{1}}^{\ast}\,. \label{eqSimilarIsomorphism}%
\end{equation}

The correspondence between $\operatorname*{End}\left(  B\left(  \mathcal{H}%
\right)  \right)  $ and representations from (\ref{eqEnd.2}) is not quite
unique: If the $\left(  S_{i}\right)  $ system is given and if $u=\left(
u_{i,j}\right)  $ is a unitary $N$-by-$N$ matrix, then the system $S_{i}%
^{u}:=\sum_{j}u_{i,j}S_{j}$ defines the same endomorphism $\alpha\left(
T\right)  =\sum_{i}S_{i}^{u}TS_{i}^{u\,\ast}$, but this is the extent of the
non-uniqueness in the correspondence.

The following result shows that the case when the induced measure
\begin{equation}
\mu_{f}\left(  \,\cdot\,\right)  =\left\Vert E\left(  \,\cdot\,\right)
f\right\Vert ^{2} \label{eqInducedMeasure}%
\end{equation}
is a product measure on
\begin{equation}
X_{N}=\Gamma_{\!N}^{\mathbb{Z}_{+}} \label{eqEnd.6}%
\end{equation}
is exceptional. Here $\mu_{f}$ is the measure defined in (\ref{eqInt.15}), and
$E\left(  \,\cdot\,\right)  $ is the projection-valued measure defined on the
Borel sets in $X_{N}$ which is induced from some given representation $\left(
S_{i}\right)  _{0\leq i<N}$ of $\mathcal{O}_{N}$. The result shows that
$\mu_{f}$ is a product measure precisely when the vector $f\in\mathcal{H}$,
$\left\Vert f\right\Vert =1$, is a simultaneous eigenvector for the operators
$S_{i}^{\ast}$. The operators $S_{i}^{\ast}$ have the form $\bigcirc
\llap{\small$\downarrow\mkern5.25mu$}$[filter], the two operator on the left
in Figure \ref{FigPerfectReconstruction}, which is the case $N=2$, i.e.,
filter followed by down-sampling; see Figure \ref{FigPerfectReconstruction} in
Section \ref{Int}.

Motivated by the main result in \cite{Jor01b}, it is appropriate to restrict
attention to irreducible representations when considering the
representations of the Cuntz algebras induced by subband filters from
signal processing: this is needed in fact for the implication 
(\ref{ProEnd.pound(1)})$\Rightarrow$(\ref{ProEnd.pound(2)})
in the proposition below. It is needed again for the states $\omega_{f}$ from
(\ref{eqEnd.10}); see also \cite{BJO04}.

\begin{prop}
\label{ProEnd.pound}\textup{(}see \textup{\cite{BJP96})} Let $\left(
S_{i}\right)  _{0\leq i<N}$ be a representation of $\mathcal{O}_{N}$ on a
Hilbert space $\mathcal{H}$, and let $\alpha=\alpha_{S}$ be the corresponding
endomorphism of $B\left(  \mathcal{H}\right)  $ \textup{(}see
\textup{(\ref{eqEnd.2}))}. Let $f\in\mathcal{H}$, $\left\Vert f\right\Vert
=1$, and let $\omega_{f}\left(  \,\cdot\,\right)  =\left\langle \, f\mid
\cdot\,f\right\rangle $ be the corresponding state. The following three
conditions are equivalent.

\begin{enumerate}
\item \label{ProEnd.pound(1)}$\omega_{f}\left(  \alpha\left(  T\right)
\right)  =\omega_{f}\left(  T\right)  $ for all $T\in B\left(  \mathcal{H}%
\right)  $.

\item \label{ProEnd.pound(2)}$f$ is a joint eigenvector for $S_{i}^{\ast}$,
$0\leq i<N$.

\item \label{ProEnd.pound(3)} There are $\lambda_{i}\in\mathbb{C}$, $0\leq
i<N$, with $\sum_{i=0}^{N} \left|  \lambda_{i} \right|  ^{2} =1$, such that%
\begin{equation}
\omega_{f}\left(  S_{i_{1}}S_{i_{2}}\cdots S_{i_{k}}S_{j_{l}}^{\ast}\cdots
S_{j_{2}}^{\ast}S_{j_{1}}^{\ast}\right)  =\bar{\lambda}_{i_{1}}\bar{\lambda
}_{i_{2}}\cdots\bar{\lambda}_{i_{k}}\lambda_{j_{1}}\lambda_{j_{2}}%
\cdots\lambda_{j_{l}} \label{eqEnd.7}%
\end{equation}
for all $k,l\in\mathbb{Z}_{+}$ and all $i_{i},\dots,i_{k}\in\Gamma_{N}$ and
$j_{1},\dots,j_{l}\in\Gamma_{N}$.\setcounter{savecounter}{\value{enumi}}
\end{enumerate}

If the representation is assumed to be type $\mathrm{I}$, then a fourth
condition is equivalent to the first three:

\begin{enumerate}
\setcounter{enumi}{\value{savecounter}}

\item \label{ProEnd.pound(4)}The measure $\mu_{f}$ obtained from $\omega_{f}$
by restriction to the maximal abelian subalgebra $\mathcal{C}$ is a product
measure on $X_{N}=\Gamma_{\!N}^{\mathbb{Z}_{+}}$.
\end{enumerate}

In general \textup{(\ref{ProEnd.pound(3)})}$\Rightarrow$%
\textup{(\ref{ProEnd.pound(4)})}.
\end{prop}

\begin{rem}
\label{RemEnd.0}If (\ref{ProEnd.pound(1)})--(\ref{ProEnd.pound(3)}) hold, the
$\lambda_{i}$ in (\ref{ProEnd.pound(3)}) are the eigenvalues from
(\ref{ProEnd.pound(2)}).
\end{rem}

\begin{rem}
\label{RemEnd.1}For non-type-$\mathrm{I}$ representations, it is possible to
have (\ref{ProEnd.pound(4)}) satisfied, even if (\ref{ProEnd.pound(1)}%
)--(\ref{ProEnd.pound(3)}) fail to hold.
\end{rem}

\begin{pf}
The equivalence of conditions (\ref{ProEnd.pound(1)}), (\ref{ProEnd.pound(2)}%
), and (\ref{ProEnd.pound(3)}) was already established in \cite{BJP96}.
Indeed, if $f$ in $\mathcal{H}$ satisfies (\ref{ProEnd.pound(2)}), there are
$\lambda_{i}\in\mathbb{C}$ such that%
\begin{equation}
S_{i}^{\ast}f=\lambda_{i}f. \label{eqEnd.8}%
\end{equation}
Using%
\[
\omega_{f}\left(  S_{i_{1}}\cdots S_{i_{k}}S_{j_{l}}^{\ast}\cdots S_{j_{1}%
}^{\ast}\right)  =\left\langle \, S_{i_{k}}^{\ast}\cdots S_{i_{1}}^{\ast
}f\bigm| S_{j_{l}}^{\ast}\cdots S_{j_{1}}^{\ast}f\right\rangle ,
\]
formula (\ref{eqEnd.7}) in (\ref{ProEnd.pound(3)}) follows. The proof that
(\ref{ProEnd.pound(1)})$\Rightarrow$(\ref{ProEnd.pound(2)}) relies on the fact
that $\omega_{f}\left(  \,\cdot\,\right)  $ is a pure state on $B\left(
\mathcal{H}\right)  $; see \cite{BJP96}. Now if the relations (\ref{eqEnd.8})
are substituted into%
\begin{equation}
\sum_{i}S_{i}S_{i}^{\ast}%
=\hbox{\upshape \small1\kern-3.8pt\normalsize1}_{\mathcal{H}}\,,
\label{eqEnd.9}%
\end{equation}
we get $\sum_{i}\left\vert \lambda_{i}\right\vert ^{2}=1$. Setting
$p_{i}=\left\vert \lambda_{i}\right\vert ^{2}$, we get a probability
distribution on $\Gamma_{N}$. Setting $k=l$ and $i_{1}=j_{1}$, $\dots$,
$i_{k}=j_{k}$ in (\ref{eqEnd.7}), we finally conclude that%
\begin{equation}
\omega_{f}\left(  S_{i_{1}}\cdots S_{i_{k}}S_{i_{k}}^{\ast}\cdots S_{i_{1}%
}^{\ast}\right)  =p_{i_{1}}p_{i_{2}}\cdots p_{i_{k}}\,. \label{eqEnd.10}%
\end{equation}
This implies that $\mu_{f}=\omega_{f}|_{\mathcal{C}}$ is a product measure.
Indeed, setting $\mu_{f}=\omega_{f}|_{\mathcal{C}}$, and using Lemma
\ref{LemPro.3}, we get
\begin{equation}
\mu_{f}\left(  \left\{  \, x\in X_{N}\mid x_{1}=i_{1},\dots,x_{k}%
=i_{k}\,\right\}  \right)  =p_{i_{1}}p_{i_{2}}\cdots p_{i_{k}}\,,
\label{eqEnd.11}%
\end{equation}
which is to say that $\mu_{f}$ is the product measure determined by the
probability distribution $\left(  p_{i}\right)  _{0\leq i<N}$. Equivalently,
$\mu_{f}$ is the unique probability measure on $X_{N}$ which satisfies the
identity (\ref{eqPro.24}).

The conclusion in (\ref{ProEnd.pound(4)}) is the fact that $\mu_{f}=\omega
_{f}|_{\mathcal{C}}$ is a product measure. Suppose the probabilities are
$\left(  p_{i}\right)  _{i\in\Gamma_{N}}$. Then (\ref{eqEnd.11}) holds. But in
terms of the $\mathcal{O}_{N}$ representation, this reads as (\ref{eqEnd.10}).

Now suppose the representation is type $\mathrm{I}$. Using again that
$\omega_{f}$ is pure, we conclude that (\ref{eqEnd.8}) must hold for some
$\lambda_{i}\in\mathbb{C}$ with $\left\vert \lambda_{i}\right\vert =p_{i}$.

To show that there are non-type-$\mathrm{I}$ representations of $\mathcal{O}%
_{N}$ for which (\ref{ProEnd.pound(4)}) holds, but (\ref{ProEnd.pound(3)})
does not, it is enough to display a state $\omega$ on $\mathcal{O}_{N}$ for
which (\ref{eqEnd.10}) holds, but (\ref{eqEnd.7}) fails. Such states are known
and studied in \cite{BJO04}. They are called KMS states. We show that for
every $p_{i}>0$, such that $\sum_{i}p_{i}=1$, there is a state $\omega
=\omega_{\left(  p\right)  }$ on $\mathcal{O}_{N}$ for which $\omega\left(
S_{i_{1}}\cdots S_{i_{k}}S_{j_{l}}^{\ast}\cdots S_{j_{1}}^{\ast}\right)
=\delta_{k,l}\delta_{i_{1},j_{1}}\cdots\delta_{i_{k},j_{k}}p_{i_{1}}\cdots
p_{i_{k}}$. Note that this is consistent with (\ref{eqEnd.10}) so
$\omega|_{\mathcal{C}}$ is a product measure on $X_{N}=\Gamma_{\!N}%
^{\mathbb{Z}_{+}}$, but it is inconsistent with (\ref{eqEnd.7}). It is known
that the representation generated by $\omega$ is type $\mathrm{III}$, i.e., it
generates a type $\mathrm{III}$ von Neumann algebra. \qed
\end{pf}

\begin{exmp}
\label{ExaEnd.3}We now calculate the measure $\mu_{f}$ from the dyadic
partitions of the unit interval $\left[  0,1\right)  $ of Example
\ref{ExaSub.1} in the case of two specific representations of $\mathcal{O}%
_{2}$. Each of the representations yields a product measure on the Cartesian
product space $X_{2}=\Gamma_{2}^{\mathbb{Z}_{+}}$: the first (a) has
probability weights $p_{0}=1$, $p_{1}=0$, and the second (b) has $p_{0}%
=p_{1}=1/2$. The induced measure on $\left[  0,1\right)  $ for (a) is the
Dirac mass $\delta_{0}$ on $\left[  0,1\right)  $, i.e., for the first
representation; and it is the restricted Lebesgue measure in the second case
(b). We take $\mathcal{H}=L^{2}\left(  \mathbb{T}\right)  $ as the Hilbert
space for both representations. We introduce the notation
\begin{equation}
e_{k}\left(  z\right)  :=z^{k} \label{eqEnd.12}%
\end{equation}
for the Fourier basis on $\mathcal{H}$, i.e., $\left\{  e_{k}\right\}
_{k\in\mathbb{Z}}$ is the usual orthonormal basis of $L^{2}\left(
\mathbb{T}\right)  $ from Fourier analysis. Following the discussion of
Section \ref{Int}, we set
\begin{equation}
S_{i}f\left(  z\right)  =m_{i}\left(  z\right)  f\left(  z^{2}\right)  ,\qquad
f\in L^{2}\left(  \mathbb{T}\right)  ,\;z\in\mathbb{T},\;i=0,1,
\label{eqEnd.13}%
\end{equation}
where $m_{0}$, $m_{1}$ is a system of quadrature mirror filters. The
conditions on the two functions may be summarized in the requirement that
\begin{equation}
\text{the matrix\quad}\frac{1}{\sqrt{2}}%
\begin{pmatrix}
m_{0}\left(  z\right)  & m_{0}\left(  -z\right) \\
m_{1}\left(  z\right)  & m_{1}\left(  -z\right)
\end{pmatrix}
\text{\quad is unitary for }\mathrm{a.a.}\,z\in\mathbb{T}. \label{eqEnd.14}%
\end{equation}
The two cases are (\ref{ExaEnd.3(1)}) and (\ref{ExaEnd.3(2)}) below.\renewcommand{\theenumi}{\alph{enumi}}

\begin{enumerate}
\item \label{ExaEnd.3(1)}\textbf{A permutative representation} (see
\cite{BrJo99a}).
\[
\text{With\quad}\left\{  {\begin{aligned}
m_{0}&=e_{0},\\
m_{1}&=e_{1},
\end{aligned}}\right\}  \text{\quad we have\quad}\left\{  {\begin{aligned}
S_{0}^{\ast}e_{0}&=e_{0},\\
S_{1}^{\ast}e_{0}&=0.
\end{aligned}}\right.
\]

\item \label{ExaEnd.3(2)}\textbf{The representation of the Haar wavelet} (see
\cite{Jor03}).
\[
\text{With\quad}\left\{  {\begin{aligned}
m_{0}&=\frac{1}{\sqrt{2}}\left( e_{0}+e_{1}\right) ,\\
m_{1}&=\frac{1}{\sqrt{2}}\left( e_{0}-e_{1}\right) ,
\end{aligned}}\right\}  \text{\quad we have\quad}\left\{  {\begin{aligned}
S_{0}^{\ast}e_{0}&=\frac{1}{\sqrt{2}}e_{0},\\
S_{1}^{\ast}e_{0}&=\frac{1}{\sqrt{2}}e_{0}.
\end{aligned}}\right.
\]

\end{enumerate}

Hence in case (\ref{ExaEnd.3(1)}), the measure $\mu_{e_{0}}$ is the Dirac mass
at $x=0$ in $\left[  0,1\right)  $, or $\mu_{e_{0}}=\delta_{0}$, and in case
(\ref{ExaEnd.3(2)}), the measure $\mu_{e_{0}}$ is Lebesgue measure on $\left[
0,1\right)  $.
\end{exmp}

The representation of $\mathcal{O}_{2}$ described in (\ref{ExaEnd.3(1)}) above
is \emph{permutative} in the sense of \cite{BrJo99a}. A permutative
representation $\left(  S_{i}\right)  _{i=0,1}$ of $\mathcal{O}_{2}$ in a
Hilbert space $\mathcal{H}$ is one for which $\mathcal{H}$ has an orthonormal
basis $\left\{  e_{n}\right\}  _{n\in\mathbb{Z}}$ such that each of the
isometries $S_{i}$ maps the basis to itself, i.e., there are maps $\sigma
_{i}\colon\mathbb{Z}\rightarrow\mathbb{Z}$ such that%
\begin{equation}
S_{i}e_{n}=e_{\sigma_{i}\left(  n\right)  },\qquad i\in\left\{  0,1\right\}
,\;n\in\mathbb{Z}. \label{eqEnd.star}%
\end{equation}
For the representation given in (\ref{ExaEnd.3(1)}), the two maps $\sigma_{i}$
are $\sigma_{0}n=2n$, $\sigma_{1}n=2n+1$. For permutative representations, the
problem of diagonalizing the commutative family of operators%
\begin{equation}
S_{i_{1}}S_{i_{2}}\cdots S_{i_{k}}S_{i_{k}}^{\ast}\cdots S_{i_{2}}^{\ast
}S_{i_{1}}^{\ast} \label{eqEnd.starstar}%
\end{equation}
is very simple; see \cite{BrJo99a}. But unfortunately wavelet representations
are typically not permutative, and the reader is referred to \cite{Jor01b} for
details of the argument.

In general, the explicit transform which diagonalizes the commuting family
(\ref{eqEnd.starstar}) may be somewhat complicated. But if the representation
$\left(  S_{i}\right)  $ is permutative, it is easy to see that the operator
monomials from (\ref{eqEnd.starstar}) may be naturally realized as
multiplication operators on the sequence space $\ell^{2}\left(  \mathbb{Z}%
\right)  $. Specifically, if $k\in\mathbb{Z}_{+}$, and $\left(  i_{1}%
,\dots,i_{k}\right)  \in\Gamma_{2}^{k}$ are given, then the corresponding
operator in (\ref{eqEnd.starstar}) is represented as multiplication by the
indicator function of the set $\sigma_{i_{1}}\sigma_{i_{2}}\cdots\sigma
_{i_{k}}\left(  \mathbb{Z}\right)  \subset\mathbb{Z}$, where the maps
$\sigma_{i}$ are determined from the formula (\ref{eqEnd.star}).

\section{\label{Com}Computation of $\mu_{f}$}

We now turn to the computation of the measures $\mu_{f}\left(  \,\cdot
\,\right)  =\left\Vert E\left(  \,\cdot\,\right)  f\right\Vert ^{2}$ in the
special case when the representation of $\mathcal{O}_{N}$ arises from a system
of subband filters. Recall from Section \ref{Pro} that every representation of
$\mathcal{O}_{N}$ defines a projection-valued measure on $\left[  0,1\right)
$ when restricted to the subalgebra $\mathcal{C}$ in $\mathcal{O}_{N}$. A
system of subband filters corresponding to $N$ subbands is a set of
$L^{\infty}$-functions $m_{0}$, $m_{1}$, $\dots$, $m_{N-1}$ on $\mathbb{T}$
such that the following matrix function on $\mathbb{T}$ takes unitary values:%
\begin{equation}
\frac{1}{\sqrt{N}}\left(  m_{j}\left(  ze^{i2\pi\frac{k}{n}}\right)
\vphantom{\frac{1}{\sqrt{N}}}\right)  _{0\leq j,k<N} \,. \label{eqCom.1}%
\end{equation}
Specifically, for $\mathrm{a.e.}\,z\in\mathbb{T}$, the $N\times N$ matrix of
(\ref{eqCom.1}) is assumed unitary.

The following lemma is well known; see \cite{BrJo02b}.

\begin{lem}
\label{LemCom.1}Let $m_{0}$, $m_{1}$, $\dots$, $m_{N-1}$ be in $L^{\infty
}\left(  \mathbb{T}\right)  $ and set%
\begin{equation}
S_{j}f\left(  z\right)  =m_{j}\left(  z\right)  f\left(  z^{N}\right)  ,\quad
f\in L^{2}\left(  \mathbb{T}\right)  ,\;z\in\mathbb{T},\;j=0,1,\dots,N-1.
\label{eqCom.2}%
\end{equation}
Then $\left(  S_{j}\right)  _{0\leq j<N}$ is a representation of
$\mathcal{O}_{N}$ on the Hilbert space $L^{2}\left(  \mathbb{T}\right)  $ if
and only if the functions $m_{j}$ satisfy the unitarity property
\textup{(\ref{eqCom.1})}.
\end{lem}

We state the next result for the middle-third Cantor set, but it applies
\emph{mutatis mutandis} to most of the fractals based on iterated function
systems (IFS's) built on affine maps.

\begin{prop}
\label{ProCom.pound}\textup{(Example \ref{ExaSub.2} revisited.)}

Let%
\begin{equation}
\left\{
{\begin{aligned} m_{0}&=\frac{1}{\sqrt{2}}\left( e_{0}+e_{2}\right) ,\\ m_{1}&=e_{1},\\ m_{2}&=\frac{1}{\sqrt{2}}\left( e_{0}-e_{2}\right) , \end{aligned}}%
\right.  \label{eqComNew.1}%
\end{equation}
and let%
\begin{equation}
S_{j}f\left(  z\right)  =m_{j}\left(  z\right)  f\left(  z^{3}\right)  ,\qquad
j=0,1,2,\;f\in L^{2}\left(  \mathbb{T}\right)  ,\;z\in\mathbb{T},
\label{eqComNew.2}%
\end{equation}
where $e_{p}\left(  z\right)  =z^{p}$, $p\in\mathbb{Z}$, and $L^{2}\left(
\mathbb{T}\right)  $ is the Hilbert space $\mathcal{H}$ of $L^{2}$-functions
on $\mathbb{T}$ defined from the Haar measure on $\mathbb{T}$. Then $\left(
S_{j}\right)  _{j=0,1,2}$ is a representation of $\mathcal{O}_{3}$ on
$\mathcal{H}=L^{2}\left(  \mathbb{T}\right)  $. Let $I$ be the unit interval,
and let $E\left(  \,\cdot\,\right)  $ be the corresponding projection-valued
measure on $\mathcal{B}\left(  I\right)  $. Then the induced scalar measure
$\mu_{e_{0}}\left(  \,\cdot\,\right)  =\left\Vert E\left(  \,\cdot\,\right)
e_{0}\right\Vert ^{2}$ is the middle-third Cantor measure of Example
\textup{\ref{ExaSub.2}}, i.e., the unique measure $\mu$ on $I$ which satisfies
\textup{(\ref{eqSub.7})}. \textup{(}It is supported on the middle-third Cantor
set $X$.\textup{)}
\end{prop}

\begin{pf}
That the system $\left(  S_{j}\right)  _{j=0,1,2}$ of (\ref{eqComNew.2}) forms
a representation of $\mathcal{O}_{3}$ on $\mathcal{H}$ follows immediately
from Lemma \ref{LemCom.1}. As noted in Lemma \ref{LemPro.2}, the corresponding
projection-valued measure $E$ is determined as follows: If $k\in\mathbb{Z}%
_{+}$, and $a=\left(  a_{1},a_{2},\dots,a_{k}
\right)  \in\Gamma_{3}^{k}$, then recall that%
\begin{equation}
E\left(  J_{k}\left(  a\right)  \right)  =S_{a}S_{a}^{\ast}%
\,,\label{eqComNew.4}%
\end{equation}
where $\displaystyle J_{k}\left(  a\right)  =\left[  \frac{a_{1}}{3}%
+\dots+\frac{a_{k}}{3^{k}},\frac{a_{1}}{3}+\dots+\frac{a_{k}}{3^{k}}+\frac
{1}{3^{k}}\right)  $, and $S_{a}:=S_{a_{1}}S_{a_{2}}\cdots S_{a_{k}}$. {}From
(\ref{eqComNew.1}), we get $\left\{  {\begin{aligned}
S_{0}^{\ast}e_{0}&=\frac{1}{\sqrt{2}}e_{0},\\
S_{1}^{\ast}e_{0}&=0,\\
S_{2}^{\ast}e_{0}&=\frac{1}{\sqrt{2}}e_{0},
\end{aligned}}\right\}  _{\mathstrut}^{\mathstrut}$, which are the joint
eigenvalue identities of (\ref{ProEnd.pound(2)}) in Proposition
\ref{ProEnd.pound}. Now a direct check on $\mu_{e_{0}}\left(  \,\cdot
\,\right)  =\left\Vert E\left(  \,\cdot\,\right)  e_{0}\right\Vert ^{2}$,
using (\ref{eqComNew.4}) and Proposition \ref{ProEnd.pound},
(\ref{ProEnd.pound(2)})$\Rightarrow$(\ref{ProEnd.pound(4)}), shows that
$\mu_{e_{0}}$ is indeed the Cantor measure of Example \ref{ExaSub.2}.
See also \cite{BJO04}. \qed
\end{pf}

For more about the representation (\ref{eqComNew.2}) and the corresponding
fractal wavelet, the reader is referred to \cite{DuJo03}. While this
representation does not correspond to a system of wavelet functions $\varphi$,
$\psi_{1}$, $\psi_{2}$ in $L^{2}\left(  \mathbb{R}\right)  $, we show in
\cite{DuJo03} that there is a Hilbert space of functions on $\mathbb{R}$ which
admits $\varphi$, $\psi_{1}$, $\psi_{2}$ as wavelet generators. If $s=\log
_{2}\left(  3\right)  =\ln2/\ln3$, the wavelet system is constructed on the
Hausdorff measure $\mathcal{H}_{s}$, i.e., the measure on $\mathbb{R}$
constructed from $\left(  dx\right)  ^{s}$ by the usual completion; see also
\cite{Fal85} for details on the Hausdorff measure $\mathcal{H}^{s}$.

\textbf{Some terminology:} For functions $g$ on $\mathbb{T}$, we define the
Fourier transform $\hat{g}\left(  n\right)  $ as follows:%
\begin{equation}
\hat{g}\left(  n\right)  =\left\langle \,e_{n}\mid g\,\right\rangle
=\int_{\mathbb{T}}z^{-n}g\left(  z\right)  \,d\lambda\left(  z\right)
=\int_{0}^{1}e^{-i2\pi n\theta}g\left(  \theta\right)  \,d\theta,
\label{eqCom.3}%
\end{equation}
where $\lambda$ denotes Haar measure on $\mathbb{T}$, and where we have
identified $g\left(  \theta\right)  $ with $g\left(  e^{i2\pi\theta}\right)  $.

If $k\in\mathbb{Z}_{+}$, and $a=\left(  a_{1},\dots,a_{k}\right)  \in
\Gamma_{\!N}^{k}$, set
\begin{equation}
m_{a}\left(  z\right)  :=m_{a_{1}}\left(  z\right)  m_{a_{2}}\left(
z^{N}\right)  \cdots m_{a_{k}}\left(  z^{N^{k-1}}\right)  , \label{eqCom.4}%
\end{equation}
or in additive notation,%
\begin{equation}
m_{a}\left(  \theta\right)  :=m_{a_{1}}\left(  \theta\right)  m_{a_{2}}\left(
N\theta\right)  \cdots m_{a_{k}}\left(  N^{k-1}\theta\right)  .
\label{eqCom.5}%
\end{equation}
When a system $m_{j}$ satisfies the condition (\ref{eqCom.1}) we say that the
representation (\ref{eqCom.2}) is a wavelet representation of $\mathcal{O}%
_{N}$.

\begin{prop}
\label{ProCom.2}Let the functions $\left(  m_{j}\right)  _{0\leq j<N}$ satisfy
the condition \textup{(\ref{eqCom.1})} and let $S_{j}$ be the corresponding
wavelet representation of $\mathcal{O}_{N}$ on the Hilbert space $L^{2}\left(
\mathbb{T}\right)  $. Let $f\in L^{2}\left(  \mathbb{T}\right)  $, $\left\Vert
f\right\Vert =1$, and let $k\in\mathbb{Z}_{+}$, $a=\left(  a_{1},\dots
,a_{k}\right)  \in\Gamma_{\!N}^{k}$. Then%
\begin{equation}
\mu_{f}\left(  J_{k}\left(  a\right)  \right)  =\sum_{n\in\mathbb{Z}%
}\left\vert \left(  f\bar{m}_{a}\right)  \widehat{\vphantom{f\bar{m}_{a}}}%
\left(  nN^{k}\right)  \right\vert ^{2}. \label{eqCom.6}%
\end{equation}

\end{prop}

\begin{pf}
Let the conditions be as stated. Then%
\begin{align*}
\mu_{f}\left(  J_{k}\left(  a\right)  \right)   &  =\left\Vert S_{a}%
S_{a}^{\ast}f\right\Vert ^{2}=\left\langle \,f\mid S_{a}S_{a}^{\ast
}f\,\right\rangle =\left\Vert S_{a}^{\ast}f\right\Vert ^{2}=\sum
_{n\in\mathbb{Z}}\left\vert \left\langle \,e_{n}\mid S_{a}^{\ast
}f\,\right\rangle \right\vert ^{2}\\
&  =\sum_{n\in\mathbb{Z}}\left\vert \left\langle \,S_{a}e_{n}\mid
f\,\right\rangle \right\vert ^{2}=\sum_{n\in\mathbb{Z}}\left\vert \left\langle
\,m_{a}\left(  z\right)  e_{n}\left(  z^{N^{k}}\right)  \bigm|f\,\right\rangle
\right\vert ^{2}\\
&  =\sum_{n\in\mathbb{Z}}\left\vert \int_{\mathbb{T}}e_{-nN^{k}}\,\bar{m}%
_{a}\,f\,d\lambda\right\vert ^{2}=\sum_{n\in\mathbb{Z}}\left\vert \left(
\bar{m}_{a}f\right)  \widehat{\vphantom{f\bar{m}_{a}}}\left(  nN^{k}\right)
\right\vert ^{2},
\end{align*}
which is the desired conclusion. \qed
\end{pf}

Specializing to $f=e_{p}$, for some $p\in\mathbb{Z}$, we get for $\mu
_{p}\left(  \,\cdot\,\right)  =\mu_{e_{p}}\left(  \,\cdot\,\right)
=\left\Vert E\left(  \,\cdot\,\right)  e_{p}\right\Vert ^{2}$,%
\begin{equation}
\mu_{p}\left(  J_{k}\left(  a\right)  \right)  =\sum_{n\in\mathbb{Z}%
}\left\vert \hat{m}_{a}\left(  p-nN^{k}\right)  \right\vert ^{2}.
\label{eqCom.7}%
\end{equation}

Let $N\in\mathbb{Z}_{+}$, $N\geq2$, and let $\left(  m_{j}\right)  _{0\leq
j<N}$ be a subband filter system, i.e., the $m_{j}$'s are functions which
satisfy condition (\ref{eqCom.1}). We shall assume further that $m_{0}$ is
Lipschitz of order $1$ as a function on $\mathbb{T}$, and that $m_{0}\left(
1\right)  =\sqrt{N}$. In that case, there is an $L^{2}\left(  \mathbb{R}%
\right)  $ scaling function $\varphi$ such that%
\[
\hat{\varphi}\left(  \xi\right)  =\prod_{k=1}^{\infty}\frac{m_{0}\left(
\xi/N^{k}\right)  }{\sqrt{N}},
\]
where we set $m_{0}\left(  \theta\right)  =m_{0}\left(  e^{-i2\pi\theta
}\right)  $, and $\hat{\varphi}\left(  \xi\right)  =\int_{\mathbb{R}}%
e^{-i2\pi\xi x}\,\varphi\left(  x\right)  \,dx$. We will assume in addition
that%
\begin{equation}
\sum_{l\in\mathbb{Z}}\left\vert \hat{\varphi}\left(  \xi+l\right)  \right\vert
^{2}=1. \label{eqCom.8}%
\end{equation}
It is known that (\ref{eqCom.8}) is equivalent to each of the following three
conditions on $\varphi$:

\begin{enumerate}
\item \label{eqCom.8equiv(1)}$\displaystyle\left\Vert \varphi\right\Vert
^{2}=\int_{\mathbb{R}}\left\vert \varphi\left(  x\right)  \right\vert
^{2}\,dx=1$,

\item \label{eqCom.8equiv(2)}the system $\displaystyle\left\{  \,\varphi
\left(  x-k\right)  \mid k\in\mathbb{Z}\,\right\}  $ is orthogonal in
$L^{2}\left(  \mathbb{R}\right)  $, and

\item \label{eqCom.8equiv(3)}the operator $\displaystyle W_{\varphi}$ is isometric;
\end{enumerate}

\noindent see \cite{BrJo02b} for details.

Recalling the $N$-adic representation for $\mathbb{Z}_{+}$,
\begin{equation}
n=a_{1}+a_{2}N+\dots+a_{k}N^{k-1},\quad k=1,2,\dots,\;a=\left(  a_{1}%
,\dots,a_{k}\right)  \in\Gamma_{\!N}^{k}, \label{eqCom.11}%
\end{equation}
we get a sequence of $L^{2}\left(  \mathbb{R}\right)  $-functions $w_{n}$, the
wavelet packet functions of Wickerhauser, satisfying%
\begin{equation}
\hat{w}_{n}\left(  \xi\right)  =N^{-k/2}m_{a_{1}}\left(  \frac{\xi}{N}\right)
m_{a_{2}}\left(  \frac{\xi}{N^{2}}\right)  \cdots m_{a_{k}}\left(  \frac{\xi
}{N^{k}}\right)  \hat{\varphi}\left(  \frac{\xi}{N^{k}}\right)  .
\label{eqCom.12}%
\end{equation}
Hence, setting $\bar{n}:=a_{k}+a_{k-1}N+\dots+a_{1}N^{k-1}$, bit-reversal, we
get
\begin{equation}
\hat{w}_{\bar{n}}\left(  N^{k}\xi\right)  =N^{-k/2}m_{a}\left(  \xi\right)
\hat{\varphi}\left(  \xi\right)  . \label{eqCom.13}%
\end{equation}
In the next lemma we shall need the following transformation $T_{\varphi}^{k}$
acting on $L^{2}\left(  \mathbb{R}\right)  $:%
\[
\left(  T_{\varphi}^{k}f\right)  \left(  x\right)  =\int_{\mathbb{R}}f\left(
x+y\right)  \,\overline{\varphi\left(  N^{k}y\right)  }\,dy.
\]

\begin{lem}
\label{LemCom.3}Let $m_{0}$, $\dots$, $m_{N-1}$ and $\varphi$ be as described
above. Then%
\begin{equation}
\hat{m}_{a}\left(  p-jN^{k}\right)  =N^{k/2}\left(  T_{\varphi}^{k}w_{\bar{n}%
}\right)  \left(  j-\frac{p}{N^{k}}\right)  . \label{eqCom.14}%
\end{equation}

\end{lem}

\begin{pf}%
\begin{align*}
T_{\varphi}^{k}w_{\bar{n}}\left(  j-\frac{p}{N^{k}}\right)   &  =\int
_{\mathbb{R}}e^{i2\pi\left(  j-\left(  p/N^{k}\right)  \right)  \xi}\left(
T_{\varphi}^{k}w_{\bar{n}}\right)  \widehat{\vphantom{f\bar{m}_{a}}}\left(
\xi\right)  \,d\xi\\
&  =N^{k}\int_{\mathbb{R}}e^{i2\pi\left(  jN^{k}-p\right)  \xi}\left(
T_{\varphi}^{k}w_{\bar{n}}\right)  \widehat{\vphantom{f\bar{m}_{a}}}\left(
N^{k}\xi\right)  \,d\xi\\
&  =\int_{\mathbb{R}}e^{i2\pi\left(  jN^{k}-p\right)  \xi}\hat{w}_{\bar{n}%
}\left(  N^{k}\xi\right)  \,\overline{\hat{\varphi}\left(  \xi\right)  }%
\,d\xi\\
&  =N^{-k/2}\int_{\mathbb{R}}e^{i2\pi\left(  jN^{k}-p\right)  \xi}m_{a}\left(
\xi\right)  \left\vert \hat{\varphi}\left(  \xi\right)  \right\vert ^{2}%
\,d\xi\\
&  =N^{-k/2}\int_{0}^{1}e^{i2\pi\left(  jN^{k}-p\right)  \xi}m_{a}\left(
\xi\right)  \sum_{l\in\mathbb{Z}}\left\vert \hat{\varphi}\left(  \xi+l\right)
\right\vert ^{2}\,d\xi\\
&  =N^{-k/2}\int_{0}^{1}e^{i2\pi\left(  jN^{k}-p\right)  \xi}m_{a}\left(
\xi\right)  \,d\xi\\
&  =N^{-k/2}\hat{m}_{a}\left(  p-jN^{k}\right)  .
\end{align*}
\qed
\end{pf}

\begin{cor}
\label{CorCom.4}Let $m_{0}$, $\dots$, $m_{N-1}$ and $\varphi$ be as described
above. Then the measure $\mu_{p}\left(  \,\cdot\,\right)  =\mu_{e_{p}}\left(
\,\cdot\,\right)  =\left\Vert E\left(  \,\cdot\,\right)  e_{p}\right\Vert
^{2}$ is given by the formula%
\begin{equation}
\mu_{p}\left(  J_{k}\left(  a\right)  \right)  =N^{k}\sum_{j\in\mathbb{Z}%
}\left\vert \left(  T_{\varphi}^{k}w_{\bar{n}}\right)  \left(  j-\frac
{p}{N^{k}}\right)  \right\vert ^{2} \label{eqCom.15}%
\end{equation}
for all $p\in\mathbb{Z}$, $k\in\mathbb{Z}_{+}$, and $a\in\Gamma_{\!N}^{k}\,$.
\end{cor}

\begin{pf}
The conclusion is immediate from the two previous lemmas, and the results in
Section \ref{End}. \qed
\end{pf}

\textbf{Some consequences of the formula (\ref{eqCom.15}):}

\begin{enumerate}
\item \label{eqCom.15conseq(1)}It gives a formula for the measure $\mu_{p}$ in
terms of the wavelet packet functions $\left(  w_{n}\right)  $ themselves. It
is known that the functions%
\begin{equation}
\left\{  N^{q/2}w_{n}\left(  N^{q}x-k\right)  \right\}  \label{eqCom.16}%
\end{equation}
form an orthonormal basis (ONB) for $L^{2}\left(  \mathbb{R}\right)  $ when
the index labels $n$, $q$, and $k$ are carefully selected: for $\left(
n,q\right)  \in\mathbb{N}\times\mathbb{Z}$ we may set $I\left(  n,q\right)
=\left[  N^{q}n,N^{q}\left(  n+1\right)  \right)  $. It is known \cite{CMW92b}
that if a subset $E$ of $\mathbb{N}\times\mathbb{Z}$ has the property that
$\left\{  \,I\left(  n,q\right)  \mid\left(  n,q\right)  \in E\,\right\}  $ is
a partition of $\left[  0,\infty\right)  $ with overlap on at most a countable
set, then
\begin{equation}
\left\{  \,N^{q/2}w_{n}\left(  N^{q}x-k\right)  \bigm|\left(  n,q\right)  \in
E,\;k\in\mathbb{Z}\,\right\}  \label{eqCom.17}%
\end{equation}
is an orthonormal basis for $L^{2}\left(  \mathbb{R}\right)  $. It is of
interest to know when the exceptional set with overlap might be more than
countable, for example if the ONB conclusion for (\ref{eqCom.17}) might hold
if it is only known that the overlap of the partition sets $\left\{
\,I\left(  n,q\right)  \mid\left(  n,q\right)  \in E\,\right\}  $ is at most
of Lebesgue measure zero: hence the interest in when the spectral measure is
absolutely continuous with respect to the Lebesgue measure on $\left[
0,1\right)  $.

\item \label{eqCom.15conseq(2)}Formula (\ref{eqCom.15}) shows that%
\[
\mu_{p+N^{k}}\left(  J_{k}\left(  a\right)  \right)  =\mu_{p}\left(
J_{k}\left(  a\right)  \right)
\]
and%
\[
\sum_{p=0}^{N^{k}-1}\mu_{p}\left(  J_{k}\left(  a\right)  \right)  =1.
\]

\item \label{eqCom.15conseq(3)}Formula (\ref{eqCom.15}) shows that the size
estimate of $\mu_{p}$ on $N$-adic intervals depends on the asymptotics of the
sequence $\left\{  \,w_{n}\mid n\in\mathbb{N}\,\right\}  $ as $n\rightarrow
\infty$, and there are effective estimates on $\left\Vert w_{n}\right\Vert
_{L^{\infty}\left(  \mathbb{R}\right)  }$ in the literature: see, e.g.,
\cite{CMW92b}, \cite{Wic93}, and \cite{Wic94}.

\item \label{eqCom.15conseq(4)}Finally, (\ref{eqCom.15}) specializes to a
known formula in case $N=2$ and $w_{n}$ is the Lemari\'{e}-Meyer wavelet
packet; see \cite{Sal03}.
\end{enumerate}

\section{\label{Fam}The family of measures $\left\{  \,\mu_{f}\mid
f\in\mathcal{H}\,\right\}  $}

Since the standard operations that are usually applied to systems of subband
filters $\left(  m_{i}\right)  _{0\leq i<N}$ depend on the functions $m_{i}$
having some degree of regularity, it is not surprising that new and different
geometric tools are needed for the analysis when the $m_{i}$'s are only known
to be measurable. In addition to the present results, the reader is referred
to recent papers of R.~Gundy, \cite{DGH00} and \cite{Gun00}.

We saw that every representation of the $C^{\ast}$-algebra $\mathcal{O}_{N}$
on a Hilbert space $\mathcal{H}$ naturally induces a family of measures
$\left\{  \,\mu_{f}\mid f\in\mathcal{H}\,\right\}  $ with each $\mu_{f}$ being
a Borel measure on the unit interval $J=\left[  0,1\right)  $. We also saw
that, if $\left(  S_{i}\right)  _{0\leq i<N}$ is a Haar wavelet representation
of $\mathcal{O}_{N}$ on $\mathcal{H}=L^{2}\left(  \mathbb{T}\right)  $, then
the measure $\mu_{e_{0}}$ is the Lebesgue measure $dt$ restricted to $J$. As
in Section \ref{End}, we denote the Fourier basis for $L^{2}\left(
\mathbb{T}\right)  $ by $e_{n}\left(  z\right)  =z^{n}$, $n\in\mathbb{Z}$,
$z\in\mathbb{T}$.

\textbf{Terminology:} Let $\mathcal{C}$ be an abelian $C^{\ast}$-algebra of
operators on a Hilbert space $\mathcal{H}$, and let $f\in\mathcal{H}$. We set
$\mathcal{C}f:=\left\{  \,Cf\mid C\in\mathcal{C}\,\right\}  $. We denote the
closure of $\mathcal{C}f$ by $\left[  \,\mathcal{C}f\,\right]  $, or just
$\mathcal{H}_{f}$ when the algebra $\mathcal{C}$ is clear from the context.

A well-known fact, based on Zorn's lemma, is that there is always a family
$f_{i}\in\mathcal{H}$, $\left\Vert f_{i}\right\Vert =1$, such that%
\begin{equation}
\mathcal{H}=\sideset{}{^{\smash{\oplus}}}{\sum}\limits_{i}\mathcal{H}_{f_{i}%
}.\label{eqFam.1}%
\end{equation}
Implicit in (\ref{eqFam.1}) is the assertion that%
\begin{equation}
\mathcal{H}_{f_{i}}\perp\mathcal{H}_{f_{j}}\text{\qquad when }i\neq
j,\label{eqFam.2}%
\end{equation}
and that the closure of the spaces $\mathcal{H}_{f_{i}}$ is all of
$\mathcal{H}$.

A vector $f$ in $\mathcal{H}$ for which $\mathcal{H}_{f}=\mathcal{H}$ is
called a \emph{cyclic vector}.

We state the next result just for the case $N=2$, but it holds for any
$N\in\mathbb{N}$, $N\geq2$.

\begin{lem}
\label{LemFam.1}Let $m_{0}=\frac{1}{\sqrt{2}}\left(  e_{0}+e_{1}\right)  $,
and $m_{1}=\frac{1}{\sqrt{2}}\left(  e_{0}-e_{1}\right)  $. Let%
\begin{equation}
S_{i}f\left(  z\right)  =m_{i}\left(  z\right)  f\left(  z^{2}\right)  ,\qquad
i=0,1,\;z\in\mathbb{T},\;f\in L^{2}\left(  \mathbb{T}\right)  ,
\label{eqFam.3}%
\end{equation}
and let $\mathcal{C}$ be the $C^{\ast}$-algebra generated by the commuting
projections%
\begin{equation}
P_{a}:=S_{a}S_{a}^{\ast} \label{eqFam.4}%
\end{equation}
as $k\in\mathbb{N}$ and $a=\left(  a_{1},\dots,a_{k}\right)  \in\Gamma_{2}%
^{k}$. Then $e_{0}$ \textup{(}i.e., the constant function $1$\textup{)} is a
cyclic vector for $\mathcal{C}$ acting on the Hilbert space $\mathcal{H}%
=L^{2}\left(  \mathbb{T}\right)  $, and the measures $\left\{  \,\mu_{f}\mid
f\in\mathcal{H}\,\right\}  $ are all absolutely continuous with respect to the
Lebesgue measure restricted to the unit interval $J=\left[  0,1\right)  $.
\end{lem}

\begin{pf}
A direct calculation, using the formula%
\begin{equation}
S_{i}^{\ast}f\left(  z\right)  =\frac{1}{2}\sum_{\substack{w\in\mathbb{T}%
\\w^{2}=z}}\overline{m_{i}\left(  w\right)  }\,f\left(  w\right)
,\label{eqFam.5}%
\end{equation}
yields%
\begin{align}
S_{i}^{\ast}e_{0} &  =\frac{1}{\sqrt{2}}e_{0},\qquad i=0,1,\label{eqFam.6}\\
S_{i}^{\ast}e_{2n} &  =\frac{1}{\sqrt{2}}e_{n},\qquad i=0,1,\label{eqFam.7}\\
S_{0}^{\ast}e_{2n+1} &  =\frac{1}{\sqrt{2}}e_{n},\;S_{1}^{\ast}e_{2n+1}%
=-\frac{1}{\sqrt{2}}e_{n}.\label{eqFam.8}%
\end{align}
More generally, if%
\begin{equation}
n=i_{1}+2i_{2}+\dots+2^{k-1}i_{k}+2^{k}p,\qquad\left(  i_{1},\dots
,i_{k}\right)  \in\Gamma_{2}^{k},\;p\in\mathbb{Z},\label{eqFam.9}%
\end{equation}
then
\begin{equation}
S_{a}^{\ast}e_{n}=\pm2^{-k/2}e_{p}\label{eqFam.10}%
\end{equation}
for $a=\left(  a_{1},\dots,a_{k}\right)  \in\Gamma_{2}^{k}$ and $S_{a}^{\ast
}:=S_{a_{k}}^{\ast}\cdots S_{a_{1}}^{\ast}$. Introducing the familiar
functions%
\begin{equation}
m_{a}\left(  z\right)  :=m_{a_{1}}\left(  z\right)  m_{a_{2}}\left(
z^{2}\right)  \cdots m_{a_{k}}\left(  z^{2^{k-1}}\right)  \label{eqFam.11}%
\end{equation}
of (\ref{eqCom.4}), we see that%
\begin{equation}
S_{a}S_{a}^{\ast}e_{0}=\pm2^{-k/2}m_{a}.\label{eqFam.12}%
\end{equation}
Let $f=\sum_{n\in\mathbb{Z}}\xi_{n}e_{n}\in L^{2}\left(  \mathbb{T}\right)  $,
and suppose $\left\langle \,f\mid S_{a}S_{a}^{\ast}e_{0}\,\right\rangle =0$
for all $k\in\mathbb{N}$ and all $a\in$ $\Gamma_{2}^{k}$. Then
\begin{equation}
\left\langle \,S_{a}^{\ast}f\mid e_{0}\,\right\rangle =0\text{ for all
}a\text{;\quad or equivalently }\int_{\mathbb{T}}S_{a}^{\ast}f\,d\lambda
=0\text{ for all }a.\label{eqFam.13}%
\end{equation}
But $S_{a}^{\ast}f=\sum_{n\in\mathbb{Z}}\xi_{n}S_{a}^{\ast}e_{n}$, and using
(\ref{eqFam.13}) and (\ref{eqFam.10}), we conclude that $\xi_{n}=0$ for all
$n\in\mathbb{Z}$, and therefore $f=0$. This means that the closed span of the
vectors%
\begin{equation}
\left\{  \,S_{a}S_{a}^{\ast}e_{0}\mid k\in\mathbb{N},\;a\in\Gamma_{2}%
^{k}\,\right\}  \label{eqFam.14}%
\end{equation}
is all of $L^{2}\left(  \mathbb{T}\right)  $. Hence, for every $h\in
L^{2}\left(  \mathbb{T}\right)  $, the space $\left[  \,\mathcal{C}h\,\right]
$ is contained in $\left[  \,\mathcal{C}e_{0}\,\right]  =L^{2}\left(
\mathbb{T}\right)  $; and the absolute continuity of $\mu_{h}$ follows from
this, since we know that $\mu_{e_{0}}$ is the Lebesgue measure on the unit
interval. \qed
\end{pf}

When the lemma is combined with the next theorem, we get the following result
for the Haar wavelet representation.

\begin{prop}
\label{ProFam.2}Let $\left(  S_{i}\right)  _{i=0,1}$ be the Haar wavelet
representation of $\mathcal{O}_{2}$ acting on $L^{2}\left(  \mathbb{T}\right)
$. Then there is a unique unitary isometry%
\begin{equation}
V\colon L^{2}\left(  \left[  0,1\right)  ^{\mathstrut},\,dt\right)
\longrightarrow L^{2}\left(  \mathbb{T}\right)  \label{eqFam.15}%
\end{equation}
such that%
\begin{equation}
V\left(  \chi_{J_{k}\left(  a\right)  }\right)  =S_{a}S_{a}^{\ast}%
e_{0}\label{eqFam.16}%
\end{equation}
for all $k\in\mathbb{N}$, $a=\left(  a_{1},\dots,a_{k}\right)  \in\Gamma
_{2}^{k}$, where $\chi_{J_{k}\left(  a\right)  }$ is the indicator function
of
\begin{equation}
J_{k}\left(  a\right)  =\left[  \frac{a_{1}}{2}+\dots+\frac{a_{k}}{2^{k}%
},\frac{a_{1}}{2}+\dots+\frac{a_{k}}{2^{k}}+\frac{1}{2^{k}}\right)
.\label{eqFam.17}%
\end{equation}
In particular, the isometry $V$ of \textup{(\ref{eqFam.15})} maps \emph{onto}
the Hilbert space $L^{2}\left(  \mathbb{T}\right)  $, and
\begin{equation}
V^{\ast}S_{a}S_{a}^{\ast}V=M_{\chi_{J_{k}\left(  a\right)  }},\label{eqFam.18}%
\end{equation}
where the operator on the right-hand side in \textup{(\ref{eqFam.18})} is
multiplication by the function $\chi_{J_{k}\left(  a\right)  }$ acting on
$L^{2}\left(  \left[  0,1\right)  ^{\mathstrut},\,dt\right)  $.
\end{prop}

\begin{thm}
\label{ThmFam.3}Let $N\in\mathbb{N}$, and let $\left\{  A_{k}\left(  a\right)
\right\}  _{k\in\mathbb{N},\;a\in\Gamma_{\!N}^{k}}$ be an $N$-adic system of
partitions of a compact metric space $X$. Let $\left(  S_{i}\right)  _{0\leq
i<N}$ be a representation of $\mathcal{O}_{N}$ on a Hilbert space
$\mathcal{H}$, and let $E^{A}\left(  \,\cdot\,\right)  $ be the corresponding
projection-valued measure, as given by Lemma \textup{\ref{LemPro.1}}.\renewcommand{\theenumi}{\alph{enumi}}

\begin{enumerate}
\item \label{ThmFam.3(1)}Then there is a set $f_{1},f_{2},\dots$
\textup{(}possibly finite\/\textup{)}, $f_{i}\in\mathcal{H}$, $\left\Vert
f_{i}\right\Vert =1$, such that the measures%
\begin{equation}
\mu_{i}\left(  \,\cdot\,\right)  :=\left\Vert E^{A}\left(  \,\cdot\,\right)
f_{i}\right\Vert ^{2}\label{eqFam.19}%
\end{equation}
are mutually singular.

\item \label{ThmFam.3(2)}For each $i$, there is a unique isometry%
\begin{equation}
V_{i}\colon L^{2}\left(  X,\mu_{i}\right)  \longrightarrow\mathcal{H}%
\label{eqFam.20}%
\end{equation}
satisfying the following three conditions:%
\begin{equation}
V_{i}\chi_{A_{k}\left(  a\right)  }=S_{a}S_{a}^{\ast}f_{i}\text{\qquad for
}k\in\mathbb{N},\;a\in\Gamma_{\!N}^{k},\label{eqFam.21}%
\end{equation}%
\begin{equation}
V_{i}^{\ast}S_{a}S_{a}^{\ast}V_{i}=M_{\chi_{A_{k}\left(  a\right)  }%
},\label{eqFam.22}%
\end{equation}
and%
\begin{equation}
V_{i}\left(  L^{2}\left(  X,\mu_{i}\right)  \right)  =\mathcal{H}_{f_{i}%
}.\label{eqFam.23}%
\end{equation}

\item \label{ThmFam.3(3)}Moreover, $\mathcal{H}%
=\sideset{}{^{\smash{\oplus}}}{\sum}\limits_{i}\mathcal{H}_{f_{i}}$, where
$\mathcal{H}_{f_{i}}\colon=\left[  \,\mathcal{C}f_{i}\,\right]  $.
\end{enumerate}
\end{thm}

\begin{pf}
The vectors $f_{i}$ may be chosen such that (\ref{ThmFam.3(3)}) holds by an
application of Zorn's lemma. With this choice, it follows from \cite{Nel69}
that the corresponding measures $\mu_{i}$ in (\ref{eqFam.19}) will be mutually singular.

When $k$ is fixed, the projections $P_{k}\left(  a\right)  =S_{a}S_{a}^{\ast}$
are mutually orthogonal, with the multi-index $a$ ranging over $\Gamma
_{\!N}^{k}$, and $E^{A}\left(  A_{k}\left(  a\right)  \right)  =P_{k}\left(
a\right)  $. Now consider the functions from (\ref{eqPro.12}). We calculate%
\begin{align*}
\int_{X}\left\vert \sum_{a}C_{a}\chi_{A_{k}\left(  a\right)  }\right\vert
^{2}\,d\mu_{i} &  =\int_{X}\sum_{a}\left\vert C_{a}\right\vert ^{2}\chi
_{A_{k}\left(  a\right)  }\,d\mu_{i}\\
&  =\int_{X}\sum_{a}\left\vert C_{a}\right\vert ^{2}\mu_{i}\left(
A_{k}\left(  a\right)  \right)  \\
&  =\int_{X}\sum_{a}\left\vert C_{a}\right\vert ^{2}\left\Vert P_{k}\left(
a\right)  f_{i}\right\Vert ^{2}\\
&  =\left\Vert \sum_{a}C_{a}P_{k}\left(  a\right)  f_{i}\right\Vert ^{2}.
\end{align*}
This proves that an isometry $V_{i}$, in (\ref{eqFam.20}), is well defined.
The argument is in fact the same measure-completion process which was used in
Section \ref{Pro}. Moreover, it follows from the construction that $V_{i}$
satisfies (\ref{eqFam.21})--(\ref{eqFam.23}). \qed
\end{pf}

\begin{ack}
We are pleased to thank Dorin Dutkay and the members of the August 2004
workshop on wavelets held at the National University of Singapore for their
interest and helpful suggestions.
We thank Professor Sandra Saliani for sending us
her recent preprints on wavelet packets and measures. We especially thank Brian
Treadway for a beautiful job of typesetting, for graphics constructions, for a
number of corrections, and for some very helpful suggestions.
We thank Professor David Larson for kindly checking the final version of 
our paper.
\end{ack}

\bibliographystyle{elsart-num}
\bibliography{jorgen}

\end{document}